


\newcount\secno
\newcount\prmno
\newcount\forno
\newif\ifnotfound
\newif\iffound

\def\namedef#1{\expandafter\def\csname #1\endcsname}
\def\nameuse#1{\csname #1\endcsname}

\long\def\ifundefined#1#2#3{\expandafter\ifx\csname
  #1\endcsname\relax#2\else#3\fi}
\def\hwrite#1#2{{\let\the=0\edef\next{\write#1{#2}}\next}}

\toksdef\ta=0 \toksdef\tb=2
\long\def\leftappenditem#1\to#2{\ta={\\{#1}}\tb=\expandafter{#2}%
                                \edef#2{\the\ta\the\tb}}
\long\def\rightappenditem#1\to#2{\ta={\\{#1}}\tb=\expandafter{#2}%
                                \edef#2{\the\tb\the\ta}}

\def\lop#1\to#2{\expandafter\lopoff#1\lopoff#1#2}
\long\def\lopoff\\#1#2\lopoff#3#4{\def#4{#1}\def#3{#2}}

\def\ismember#1\of#2{\foundfalse{\let\given=#1%
    \def\\##1{\def\next{##1}%
    \ifx\next\given{\global\foundtrue}\fi}#2}}

\def\section#1{\medbreak
               \global\def\currenvir{section}
               \global\advance\secno by1\global\prmno=0
               {\bf \number\secno. {#1}}
               \smallskip}

\def\subsection{\global\def\currenvir{subsection}
                \global\advance\prmno by1
                \ind{ (\number\secno.\number\prmno) }}
\def\subsec{\global\def\currenvir{subsection}
                \global\advance\prmno by1
                { (\number\secno.\number\prmno)\ }}

\def\formule{\global\def\currenvir{formule}
                \global\advance\forno by1\global\subsecno=0
            }

\def\subformule{\global\def\currenvir{subformule}
                \global\advance\subsecno by1
            }

\def\proclaim#1{\global\advance\prmno by 1
                {\bf #1 \the\secno.\the\prmno$.-$ }}

\long\def\th#1 \enonce#2\endth{%
   \medbreak\proclaim{#1}{\it #2}\global\def\currenvir{th}\smallskip}

\def\rem#1{\global\advance\prmno by 1
{\it #1} \the\secno.\the\prmno$.-$}


\def\isinlabellist#1\of#2{\notfoundtrue%
   {\def\given{#1}%
    \def\\##1{\def\next{##1}%
    \lop\next\to\za\lop\next\to\zb%
    \ifx\za\given{\zb\global\notfoundfalse}\fi}#2}%
    \ifnotfound{\immediate\write16%
                 {Warning - [Page \the\pageno] {#1} No reference found}}%
                \fi}%
\def\ref#1{\ifx\labellist\empty{\immediate\write16
                 {Warning - No references found at all.}}
               \else{\isinlabellist{#1}\of\labellist}\fi}

\def\newlabel#1#2{\rightappenditem{\\{#1}\\{#2}}\to\labellist}
\def\labellist{}

\def\label#1{\relax}

\def\openall{\openout\lbl=\jobname.lbl}

\newread\testfile
\def\lookatfile#1{\openin\testfile=\jobname.#1
    \ifeof\testfile{\immediate\openout\nameuse{#1}\jobname.#1
                    \write\nameuse{#1}{}
                    \immediate\closeout\nameuse{#1}}\fi%
    \immediate\closein\testfile}%

\def\begin{\newlabel{comp}{1.2}
\newlabel{Bun_G}{1.3}
\newlabel{prod}{1.4}
\newlabel{Pic-Bun}{1.5}
\newlabel{twist}{2}
\newlabel{M_G^\delta }{2.1}
\newlabel{deg}{2.2}
\newlabel{ex-M_G^\delta }{2.3}
\newlabel{unif}{2.4}
\newlabel{PGL_r}{3.1}
\newlabel{Sp^d}{4.1}
\newlabel{PSp}{4.2}
\newlabel{\theta }{5.2}
\newlabel{PSO}{5.3}
\newlabel{ext}{6}
\newlabel{list}{6.2}
\newlabel{char}{6.3}
\newlabel{produit}{6.4}
\newlabel{section-mod}{6.5}
\newlabel{desc-prelim}{6.6}
\newlabel{desc-cor}{6.7}
\newlabel{desc-prop}{6.10}
\newlabel{Theta }{6.11}
\newlabel{restric}{6.12}
\newlabel{coarse}{7.1}
\newlabel{unirat}{7.2}
\newlabel{comppres}{7.3}
\newlabel{cyc}{7.4}
\newlabel{Pic(M_Spin)}{8}
\newlabel{Pic_Spin}{8.2}
\newlabel{G_2}{8.4}
\newlabel{A0}{9.1}
\newlabel{BC0}{9.3}
\newlabel{A_d}{10.1}
\newlabel{slopes}{10.2}
\newlabel{KP}{10.3}
\newlabel{desc-KP}{10.4}
\newlabel{M(2l,l)}{10.6}
\newlabel{C}{11.2}
\newlabel{e_G}{12.1}
\newlabel{factorial}{13.3}
\newlabel{aut}{13.4}}


\magnification 1250
\mathcode`A="7041 \mathcode`B="7042 \mathcode`C="7043
\mathcode`D="7044 \mathcode`E="7045 \mathcode`F="7046
\mathcode`G="7047 \mathcode`H="7048 \mathcode`I="7049
\mathcode`J="704A \mathcode`K="704B \mathcode`L="704C
\mathcode`M="704D \mathcode`N="704E \mathcode`O="704F
\mathcode`P="7050 \mathcode`Q="7051 \mathcode`R="7052
\mathcode`S="7053 \mathcode`T="7054 \mathcode`U="7055
\mathcode`V="7056 \mathcode`W="7057 \mathcode`X="7058
\mathcode`Y="7059 \mathcode`Z="705A
\def\spacedmath#1{\def\packedmath##1${\bgroup\mathsurround =0pt##1\egroup$}
\mathsurround#1
\everymath={\packedmath}\everydisplay={\mathsurround=0pt}}
\def\nospacedmath{\mathsurround=0pt
\everymath={}\everydisplay={} } \spacedmath{2pt}
\font\eightrm=cmr8         \font\eighti=cmmi8
\font\eightsy=cmsy8        \font\eightbf=cmbx8
        \font\eightit=cmti8
        \font\sixrm=cmr6
\font\sixi=cmmi6           \font\sixsy=cmsy6
\font\sixbf=cmbx6
\catcode`\@=11
\def\eightpoint{%
  \textfont0=\eightrm \scriptfont0=\sixrm \scriptscriptfont0=\fiverm
  \def\rm{\fam\z@\eightrm}%
  \textfont1=\eighti  \scriptfont1=\sixi  \scriptscriptfont1=\fivei
    \textfont2=\eightsy \scriptfont2=\sixsy \scriptscriptfont2=\fivesy
   \textfont\itfam=\eightit
  \def\it{\fam\itfam\eightit}%
   \textfont\bffam=\eightbf \scriptfont\bffam=\sixbf
  \scriptscriptfont\bffam=\fivebf
  \def\bf{\fam\bffam\eightbf}%
   \abovedisplayskip=9pt plus 3pt minus 9pt
  \belowdisplayskip=\abovedisplayskip
  \abovedisplayshortskip=0pt plus 3pt
  \belowdisplayshortskip=3pt plus 3pt
  \smallskipamount=2pt plus 1pt minus 1pt
  \medskipamount=4pt plus 2pt minus 1pt
  \bigskipamount=9pt plus 3pt minus 3pt
  \normalbaselineskip=9pt
  \setbox\strutbox=\hbox{\vrule height7pt depth2pt width0pt}%
  \let\bigf@nt=\eightrm     \let\smallf@nt=\sixrm
  \normalbaselines\rm}

\def\qfl#1{\buildrel {#1}\over {\longrightarrow}}
\def\phfl#1#2{\normalbaselines{\baselineskip=0pt
\lineskip=10truept\lineskiplimit=1truept}\nospacedmath\smash
{\kern-4pt\mathop{\hbox
to 8truemm{\rightarrowfill}\kern-4pt}
\limits^{\scriptstyle#1}_{\scriptstyle#2}}}
\def\hfl#1#2{\normalbaselines{\baselineskip=0truept
\lineskip=10truept\lineskiplimit=1truept}\nospacedmath\smash{\mathop{\hbox to
12truemm{\rightarrowfill}}\limits^{\scriptstyle#1}_{\scriptstyle#2}}}
\def\diagramme#1{\def\normalbaselines{\baselineskip=0truept
\lineskip=10truept\lineskiplimit=1truept}   \matrix{#1}}
\def\vfl#1#2{\llap{$\scriptstyle#1$}\left\downarrow\vbox to
6truemm{}\right.\rlap{$\scriptstyle#2$}}
\def\note#1#2{\footnote{\parindent
0.4cm$^#1$}{\vtop{\eightpoint\baselineskip12pt\hsize15.5truecm\noindent #2}}
\parindent 0cm}
\def\mono{\lhook\joinrel\mathrel{\longrightarrow}}
\def\iso{\mathrel{\mathop{\kern 0pt\longrightarrow }\limits^{\sim}}}

\def\sdir_#1^#2{\mathrel{\mathop{\kern0pt\oplus}\limits_{#1}^{#2}}}
\def\pprod_#1^#2{\raise
2pt
\hbox{$\mathrel{\scriptstyle\mathop{\kern0pt\prod\,}\limits_{#1}^{#2}}$}}

\def\pc#1{\tenrm#1\sevenrm}
\def\up#1{\raise 1ex\hbox{\smallf@nt#1}}
\def\tx{\kern-1.5pt -}
\def\cqfd{\kern 2truemm\unskip\penalty 500\vrule height 4pt depth 0pt width
4pt\medbreak} 
\def\virg{\raise
.4ex\hbox{,}}
\def\decale#1{\smallbreak\hskip 28pt\llap{#1}\kern 5pt}
\def\no{n\up{o}\kern 2pt}
\def\ind{\par\hskip 1truecm\relax}
\def\indp{\par\hskip 0.5cm\relax}
\def\moins{\mathrel{\hbox{\vrule height 3pt depth -2pt width 6pt}}}
\def\rond{\kern 1pt{\scriptstyle\circ}\kern 1pt}
\def\iso{\mathrel{\mathop{\kern 0pt\longrightarrow }\limits^{\sim}}}

\def\Hom{\mathop{\rm Hom}\nolimits}
\def\Aut{\mathop{\rm Aut}\nolimits}

\def\Id{\mathop{\rm Id}\nolimits}
\def\Ker{\mathop{\rm Ker}\nolimits}

\def\Spec{\mathop{\rm Spec}}
\def\det{\mathop{\rm det}\nolimits}
\def\Pic{\mathop{\rm Pic}\nolimits}
\def\Div{\mathop{\rm Div}\nolimits}
\def\dim{\mathop{\rm dim}\nolimits}

\def\div{\mathop{\rm div\,}\nolimits}

\def\ad{\mathop{\rm Ad}\nolimits}

\font\gragrec=cmmib10
\def\Mu{\hbox{\gragrec \char22}}
\font\san=cmssdc10
\def\wedge{\hbox{\san \char3}}
\def\sym{\hbox{\san \char83}}
\def\sdir_#1^#2{\mathrel{\mathop{\kern0pt\oplus}\limits_{#1}^{#2}}}

\input amssym.def
\input amssym
\vsize = 25truecm
\hsize = 16truecm
\voffset = -.5truecm
\parindent=0cm
\baselineskip15pt
\overfullrule=0pt

\let\bk\backslash
\let\lra\longrightarrow
\let\ra\rightarrow
\font\gros=cmbx12
\font\tgros=cmbx12 at 14pt
\begin
\null\vskip0.5cm
\centerline{\tgros The Picard group of the moduli of G-bundles on a
curve}

\smallskip
\centerline{Arnaud {\pc BEAUVILLE}\note{1}{Partially supported by the
European HCM project ``Algebraic Geometry in Europe" (AGE).}, Yves  {\pc
LASZLO}$^1$, Christoph {\pc SORGER}\note{2}{Partially supported by
Europroj.}}

\vskip 1cm
{\bf Introduction}

\ind This paper is concerned with the moduli space of principal $G$\tx
bundles on an algebraic curve, for $G$ a complex semi-simple group.
 While the case $G={\bf SL}_r$, which corresponds to vector bundles, has
been extensively studied in algebraic geometry, the general case has
attracted much less attention until recently, when it became clear that
these spaces play an important role  in Quantum Field Theory. In
particular, if $L$ is a holomorphic line bundle on the moduli space $M_G$,
the space
$H^0(M_G,L)$ is essentially independent of the curve $X$, and can be
naturally identified with what physicists call the {\it space of conformal
blocks} associated to the most standard  Conformal Field Theory,
the so-called  WZW-model. This gives
a strong motivation to determine the  group $\Pic(M_G)$ of holomorphic line
bundles on the moduli space.

\ind Up to this point we have been rather vague about what we
should call the moduli space of $G$\tx bundles on $X$. Unfortunately
there are two possible choices, and both are meaningful. Because $G$\tx
bundles have usually nontrivial automorphisms, the natural solution to the
moduli problem is not an algebraic variety, but a slightly more complicated
object, the algebraic stack ${\cal M}_G$. This has all the good
properties one
 expects from a moduli space; in particular, a line bundle on ${\cal M}_G$
is the functorial assignment, for every variety $S$ and every $G$\tx bundle
on $X\times S$, of a line bundle on $S$.
There is also a more down-to-earth object, the coarse moduli space
$M_G$ of semi-stable $G$\tx bundles; the group $\Pic(M_G)$ is a subgroup of
$\Pic({\cal M}_G)$, but its geometric meaning is less clear.
\ind In this paper we determine the groups  $\Pic(M_G)$ and
$\Pic({\cal M}_G)$ for essentially all classical semi-simple groups,
i.e.\ of type $A,B,C,D$ and $G_2$. Since the simply-connected case was
treated in [L-S] (see also [K-N]), we are mainly concerned with  non
simply-connected groups. One new difficulty appears: the moduli space is no
longer connected,  its connected components are naturally indexed by $\pi
_1(G)$. Let $\widetilde{G}$ be the universal covering of $G$; for each
$\delta\in \pi _1(G)$, we construct a natural ``twisted" moduli stack ${\cal
M}_{\widetilde{G}}^\delta$ which dominates ${\cal M}_G^\delta$. (For
instance if $G={\bf PGL}_r$,  it is the moduli stack of vector bundles on
$X$ of rank $r$ and fixed determinant of degree $d$, with $e^{2\pi
id/r}=\delta$.) This moduli stack carries in each case a natural  line
bundle ${\cal D}$, the determinant bundle associated to the standard
representation of $\widetilde{G}$. We can now state some of our results; for
simplicity we only consider the adjoint groups.
\medskip
{\bf Theorem}$.-$ {\it Put
$\varepsilon _G^{}=1$ if  the rank of $G$ is even, $2$ if it is odd. Let
$\delta
\in \pi _1(G)$. \indp{\rm a)} The torsion subgroup of $\Pic({\cal
M}_G^\delta)$ is isomorphic to $H^1(X,\pi _1(G))$. The torsion-free quotient
is infinite cyclic, generated by ${\cal
D}^r$ if $G={\bf PGL}_r$, by ${\cal D}^{\varepsilon_G^{}}$ if $G={\bf
PSp}_{2l}$ or ${\bf PSO}_{2l}$.} \indp b) {\it The group $\Pic(M_G^\delta)$
is infinite cyclic, generated by ${\cal D}^{r\varepsilon_G^{}}$ if $G={\bf
PGL}_r$, by ${\cal D}^{2\varepsilon_G^{}}$ if $G={\bf PSp}_{2l}$ or ${\bf
PSO}_{2l}$.}\note{1} {The statement
``$\Pic(M_G)$ is generated by ${\cal D}^k$"
 must be interpreted as ``${\cal D}^k$ descends to $M_G$, and
the line bundle on $M_G$ thus obtained generates $\Pic(M_G)$" -- and
similarly for a).}
\medskip

\ind Unfortunately, though our method has some general features, it
requires  a  case-by-case analysis -- in view of the result, this is perhaps
unavoidable. An amusing consequence (\S 13) is that the moduli space $M_G$
is {\it not} locally factorial, except when $G$ is simply connected with
each simple factor of type $A,C$ or perhaps $E$. However it is always a
Gorenstein variety.

\vskip1,5cm
{\bf Notation}\smallskip
\ind  Throughout this paper we denote by $X$  a smooth projective connected
curve over ${\bf C}$; we fix a point $p$ of $X$. We let  $G$ be a complex
semi-simple group; by a $G$\tx bundle we always mean a principal bundle
with structure group $G$.  We denote by ${\cal M}_G$ the moduli stack
parameterizing
$G$\tx  bundles  on $X$, and by $M_G$ the coarse moduli variety of
semi-stable $G$\tx bundles (see \S 7).
\vskip2cm
\centerline{\gros Part I: The Picard group of the moduli stack}

\vskip1cm

\section{The stack ${\cal M}_G$}

\subsection   Our
main tool to  study $\Pic({\cal M}_G)$  will be the uniformization theorem
of [B-L], [F2] and [L-S], which we now recall. We denote by $LG$ the loop
group $G({\bf C}((z)))$, viewed as an ind-scheme over ${\bf C}$,  by $L^+G$
the sub-group scheme
$G({\bf C}[[z]])$, and by ${\cal Q}_G$ the infinite Grassmannian $LG/L^+G$;
it is a  direct limit of projective integral varieties ({\it loc.\ cit.}).
Finally let $L_XG$ be the sub-ind-group $G({\cal O}(X\moins p))$ of $LG$.
The uniformization theorem defines a canonical isomorphism of stacks
$${\cal M}_G\ \iso\ L_XG\bk {\cal Q}_G\ .$$

\ind Let $\widetilde{G}\rightarrow G$ be the universal cover of $G$; its kernel
is canonically isomorphic to $\pi_1(G)$.  We want to compare the stacks
${\cal M}_G$ and ${\cal M}_{\widetilde{G}}$.
\th Lemma
\enonce
{\rm (i)} The group $\pi_0(LG)$  is canonically isomorphic to
$\pi_1(G)$.

{\rm (ii)}  The quotient map $LG\rightarrow {\cal Q}_G$ induces a
bijection  $\pi_0(LG)\rightarrow \pi_0({\cal Q}_G)$. Each connected component
of ${\cal
Q}_G$ is isomorphic to  ${\cal Q}_{\widetilde{G}}$.

{\rm (iii)} The group $\pi_0(L_XG)$ is canonically isomorphic to
$H^1(X,\pi_1(G))$.

{\rm (iv)} The group $L_XG$ is contained in the neutral component $(LG)^{\rm
o}$ of
$LG$.

\endth
\label{comp}
{\it Proof}: Let us first prove (i) when $G$ is simply connected. In that
case, there exists a finite family  of
homomorphisms $x_\alpha:{\bf G}_a\rightarrow G$ such that for any extension
$K$ of ${\bf C}$, the subgroups $x_\alpha(K)$ generate $G(K)$ [S1]. Since the
ind-group ${\bf G}_a({\bf C}((z)))$ is connected, it follows that $LG$ is
connected.
\ind In the general case, consider the exact sequence
$1\ra \pi_1(G)\ra \widetilde{G}\ra G\ra 1$ as an exact sequence of \'etale
sheaves on  $D^*:=\Spec {\bf C}((z))$. Since
$H^1(D^*,\widetilde{G})$ is trivial [S2], it gives rise to an exact
sequence of ${\bf C}$\tx groups
$$1\ra L\widetilde{G}/\pi_1(G)\lra LG\lra H^1(D^*,\pi_1(G))\ra
1\ .\leqno{(\ref{comp}\  {\it a})}$$
The  assertion  (i) follows from the connectedness of
$L\widetilde{G}$  and the canonical isomorphism $H^1(D^*,\pi_1(G))\iso
\pi_1(G)$ (Puiseux theorem).

\ind To prove (ii), we first observe that the group $L^+G$ is connected: for
any
$\gamma \in L^+G({\bf C})$, the map $ F_\gamma :G\times {\bf
A}^1\rightarrow L^+G $ defined by $ F_\gamma (g,t)=g^{-1}\gamma (tz) $
satisfies $F_\gamma (\gamma (0),0)= 1$ and $F_\gamma (1,1)=\gamma$, hence
connects
$\gamma $ to the origin. Therefore the canonical map $\pi_0(LG)\rightarrow
\pi_0(LG/L^+G)$ is bijective. Moreover it follows from (\ref{comp} {\it a})
that
$(LG)^{\rm o}$ is isomorphic to  $L\widetilde{G}/\pi_1(G)$, which gives (ii).

\ind  Consider now the cohomology exact
sequence on $X^*$ associated to  the exact sequence $1\ra
\pi_1(G)\ra \widetilde{G}\ra G\ra 1$. Since $H^1(X^*,\widetilde{G})$ is trivial
[Ha], we get  an exact sequence of ${\bf C}$\tx groups  $$1\ra
L_X\widetilde{G}/\pi_1(G)\ra L_XG\ra H^1(X^*,\pi_1(G))\ra 1\
.\leqno(\ref{comp}\ {\it b})$$  Since the restriction map
$H^1(X,\pi_1(G))\rightarrow H^1(X^*,\pi_1(G))$ is bijective and
$L_X\widetilde{G}$ is connected ([L-S], Prop.\ 5.1), we obtain (iii). \ind
Comparing (\ref{comp}\ {\it a}) and (\ref{comp}\ {\it b}) we see that (iv) is
equivalent to saying that the restriction map $H^1(X^*,\pi_1(G))\rightarrow
H^1(D^*,\pi_1(G))$ is zero. This follows at once from the commutative diagram
of restriction maps $$\diagramme{H^1(X,\pi_1(G)) & \hfl{\sim}{} &
H^1(X^*,\pi_1(G))\cr \vfl{}{} & & \vfl {}{} \cr H^1(D,\pi_1(G)) & \hfl{}{} &
H^1(D^*,\pi_1(G)) }$$and the vanishing of $H^1(D,\pi_1(G))$.\cqfd

\medskip
\ind For $\delta\in \pi_1(G)$, let us denote by $(LG)^\delta$ the component of
$LG$
corresponding to $\delta$ via Prop.\ \ref{comp} (i).

\th Proposition
\enonce {\rm a)} There is a canonical bijection $\pi_0({\cal M}_G)\iso
\pi_1(G)$.
\ind {\rm b)} For $\delta\in \pi_1(G)$, let ${\cal M}_G^\delta$ be the
corresponding
component of ${\cal M}_G$; let $\zeta $ be any element of $(LG)^\delta({\bf
C})$.
There is a canonical isomorphism $${\cal M}_G^\delta\ \iso\
(\zeta^{-1}\,L_XG\,\zeta)\bk {\cal Q}_{\widetilde{G}}\ .$$
\endth
\label{Bun_G}
{\it Proof}: The first assertion  follows from the uniformization theorem and
Lemma
\ref{comp}, (i), (ii) and (iv). Again by  the uniformization theorem,
${\cal M}_G^\delta$  is isomorphic to
$L_XG\bk (LG)^\delta/L^+G$; left multiplication by $\zeta^{-1}
$ induces an isomorphism of $(LG)^\delta/L^+G$ onto $(LG)^{\rm o}/L^+G={\cal
Q}_{\widetilde{G}}$, and therefore an isomorphism of $L_XG \bk
(LG)^\delta/L^+G$ onto
$(\zeta^{-1}\,L_XG\,\zeta)\bk  {\cal Q}_{\widetilde{G}}$.\cqfd

\medskip

  \ind  Prop.\ \ref{Bun_G} a)  assigns to any $G$\tx bundle $P$ on $X$ an
element
$\delta$ of $\pi_1(G)$ such that $P$ defines a point of ${\cal M}_G^\delta$; we
will
refer to $\delta$ as the {\it degree} of $P$. \ind  We will use  Prop.\
\ref{Bun_G}  to determine the Picard group of ${\cal M}_G^\delta$; therefore we
first
need  to compute $\Pic({\cal Q}_{\widetilde{G}})$. We denote by $s$ the number
of
simple factors of ${\rm Lie}(G)$.
\th  Lemma \enonce  The Picard group of ${\cal
Q}_{\widetilde{G}}$ is isomorphic to ${\bf Z}^s$.   \endth\label{prod}

{\it Proof}:  Write $\widetilde{G}$ as a product
$\pprod_{i=1}^s\widetilde{G}_i$ of almost simple simply connected groups.
Put ${\cal Q}={\cal Q}_{\widetilde{G}}$ and ${\cal Q}_i={\cal
Q}_{\widetilde{G}_i}$;  the Grassmannian ${\cal Q}$ is isomorphic to
$\pprod_{}^{}{\cal Q}_i$. The Picard group of ${\cal Q}_i$ is free of rank
$1$ [M]; we denote by ${\cal O}_{{\cal Q}_i}(1)$ its positive generator. The
projections ${\cal Q}\rightarrow {\cal Q}_i$ define a group homomorphism
$\pprod_{}^{}\Pic({\cal Q}_i)\rightarrow \Pic({\cal Q})$; we claim that it is
bijective.    \ind Let ${\cal L}$ be a line bundle on ${\cal Q}$; there are
integers $(m_i)$ such that the restriction of ${\cal L}$ to
$\{q_1\}\times\ldots\times{\cal Q}_j \times\ldots\times\{q_s\}$, for any
$(q_i)\in \pprod_{}^{}{\cal Q}_i$ and any $j$, is isomorphic to ${\cal
O}_{{\cal
Q}_j}(m_j)$.  Then ${\cal L}$ is isomorphic to  $\boxtimes_i\,{\cal
O}_{{\cal Q}_i}(m_i)$: by writing   each ${\cal Q}_i$ as a direct limit of
varieties ${\cal Q}_i^{(n)}$, we are reduced to prove that these two line
bundles are isomorphic over $\pprod_i^{}{\cal Q}_i^{(n)}$, which follows
immediately from the theorem of the square.\cqfd

\medskip
\ind If $A$ is a finite abelian group, we will denote by $\widehat{A}$ its
Pontrjagin dual\break $\Hom(A,{\bf C}^*)$; it is  isomorphic
(non-canonically) to $A$. \medskip
\th Proposition
\enonce
For $\delta\in \pi_1(G)$, let  $q_G^\delta:{\cal
Q}_{\widetilde{G}}\rightarrow {\cal M}_G^\delta$ be the canonical projection
{\rm
(Prop.\ \ref{Bun_G}).} The kernel of the homomorphism
 $$ (q_G^\delta)^*:\Pic({\cal M}_G^\delta)
\longrightarrow  \Pic({\cal Q}_{\widetilde{G}})\cong {\bf  Z}^{s}$$
is canonically isomorphic to $ H^1(X, \pi_1(G)\,\widehat{}\ )$, and its
image  has finite index. \endth\label{Pic-Bun}

 {\it Proof}:  Since $q_G^\delta$
identifies ${\cal M}_G^\delta$ to the quotient of ${\cal Q}_{\widetilde{G}}$ by
$\zeta^{-1}\,L_XG\,\zeta$, line bundles on ${\cal M}_G^\delta$ correspond in a
one-to-one way to line bundles on ${\cal Q}_{\widetilde{G}}$ with a
$(\zeta^{-1}\,L_XG\,\zeta)$\tx linearization; in particular, the kernel of
$(q_G^\delta)^*$ is
canonically isomorphic to the character group  $\Hom(L_XG,{\bf C}^*)$. From the
exact
sequence (\ref{comp} {\it b})  and the triviality of the character group of
$L_X\widetilde{G}$ ([L-S], Cor.\ 5.2) we see that the group $\Hom(L_XG,{\bf
C}^*)$ is
isomorphic to  $H^1(X, \pi_1(G))\ \widehat{}\ $, which  can
 be identified by duality with $H^1(X, \pi_1(G)\,\widehat{}\ )$.
\ind Write  $\widetilde{G}\cong\pprod_{i=1}^s\widetilde{G}_i$ as in Lemma
\ref{prod}. The image of $\pi_1(G)$ under the $i$\tx th projection
$p_i:\widetilde{G}\rightarrow \widetilde{G}_i$ is a central subgroup $A_i$
of $\widetilde{G}_i$; we denote by $G_i$ the quotient $\widetilde{G}_i/A_i$,
so that $p_i$ induces a homomorphism
  $G\rightarrow G_i$. Let $\delta_i$ be the image of $\delta$ in
$\pi_1(G_i)$.  Choosing   a non trivial representation $\rho:G_i\rightarrow
{\bf SL}_r $ gives rise to  a commutative diagram
$$\diagramme{
{\cal Q}_{\widetilde{G}} &\hfl{pr_i}{}
&{\cal Q}_{\widetilde{G}_i}&\hfl{}{}& {\cal Q}_{{\bf SL}_r} &  \cr
 \vfl{q_G^\delta}{} & & \vfl{q_{G_i}^{\delta_i}}{} &  &  \vfl{q^{}_{{\bf
SL}_r}}{}
&\cr  {\cal M}_G^\delta &\hfl{}{} &
{\cal M}_{G_i}^{\delta_i} & \hfl{}{} & {\cal M}_{{\bf SL}_r}& \kern-10pt .
}$$
The pull back of the determinant bundle ${\cal D}$ on ${\cal M}_{{\bf
SL}_r}$ to ${\cal Q}_{{\bf SL}_r}$ is ${\cal O}_{\cal Q}(1)$ [B-L], and  the
pull back of ${\cal O}_{\cal Q}(1)$ to ${\cal Q}_i:={\cal Q}_{\widetilde{G}_i}$
is
${\cal O}_{{\cal Q}_i}(d_\rho )$ for some integer $d_\rho $ (the Dynkin index
of
$\rho $, see [L-S]). Therefore $pr_i^*\,{\cal O}_{{\cal Q}_i}(d_\rho )$ belongs
to the
image of $(q_G^\delta)^*$. It follows that this image has finite index.\cqfd

\rem{Remark}  In the sequel we
will be mostly interested in the case where $G$ is  almost simple; then
 $\pi_1(G)$ is canonically isomorphic to $\Mu_n$ (the group of $n$\tx th
roots of $1$) or to $\Mu_2\times \Mu_2$, and each of these groups is
naturally isomorphic to its dual (by  choosing $e^{2\pi i/n}$ as generator
of $\Mu_n$). We thus get that the torsion subgroup of  $\Pic({\cal
M}_G^\delta)$ is $J_n$ in the first case and $J_2\times J_2$ in the second,
where $J_n$ denotes the kernel of the multiplication by $n$ in the Jacobian
of $X$.

\vskip1cm
\section {The twisted moduli stack ${\cal M}_G^\delta$}
\label{twist}
\subsection \label{M_G^\delta} Proposition \ref{Pic-Bun} takes care of the
torsion subgroup of $\Pic({\cal M}_G^\delta)$; to complete the description
of this group we need to determine the image of $(q_G^\delta)^*$, or more
precisely to describe geometrically the generators of this image. To do
this we will again compare with the simply connected case, by constructing
for every $\delta\in \pi_1(G)$ a ``twisted" moduli stack  ${\cal
M}_{\widetilde{G}}^\delta$  which dominates ${\cal M}_G^\delta$.
\ind  Let  $A$ be a central subgroup of $G$, together with an
isomorphism $ A\iso  \pprod_{j=1}^s\Mu_{r_j}$.  Using  this
isomorphism we  identify  $A$ to a subgroup of  the torus $T=({\bf G}_m)^s$;
let $C_AG $ be the quotient of $G\times T$ by the diagonal subgroup $A$.
The projection $\partial:C_AG    \rightarrow  T/A\cong T$ induces a
morphism of stacks $\det:{\cal M}_{C_AG    }\rightarrow {\cal M}_{T}$. For
each element ${\bf d}=(d_1,\ldots,d_s)$ of
${\bf Z}^s$, let us denote by ${\cal O}_X({\bf d}p)$ the rational point of
${\cal M}_T$ defined by
$({\cal O}_X(d_1p),\ldots,{\cal O}_X(d_sp))$.
The fiber ${\cal M}^{{\bf d}}_{G,A}$ of $det$ at
${\cal O}_X({\bf d}p)$ depends only, up to a canonical isomorphism,  of the
class of ${\bf d}$ modulo ${\bf r}=(r_1,\ldots,r_s)$.

\ind  If $S$
is a complex scheme,  an object of ${\cal M}^{\bf d}_{G,A}(S)$
is by definition a
$C_AG    $\tx bundle $P$ on $X\times S$ together with
 a  $T$\tx bundle isomorphism of
 $P\times^{C_AG    }T$  with  the $T$\tx bundle associated to ${\cal
O}_X({\bf d}p)$. If
${\bf d}=0$, giving such an isomorphism amounts to reduce the structure group
of $P$ to $\Ker\partial=G$: in other
words, the stack ${\cal M}^{\bf 0}_{G,A} $ is canonically isomorphic to
${\cal M}_G$.

\subsection\label{deg}
The projection $p:C_AG    \rightarrow G/A$ induces a  morphism of
stacks\break $\pi:{\cal M}_{G,A}^{\bf d}\rightarrow{\cal M}_{G/A}$.
The exact sequence
$$1\rightarrow A\longrightarrow C_AG    \ \hfl{(p,\partial)}{}\ (G/A)\times
T\rightarrow 1$$
gives rise to a cohomology exact sequence
$$H^1(X,A)\rightarrow H^1(X,C_AG    )\rightarrow H^1(X,G/A)\times
H^1(X,T)\rightarrow H^2(X,A)$$ from which we deduce that the degree
$\delta\in \pi_1(G)$ of the $\!G$\tx bundle $\pi(P)$, for  $P\in$ ${\cal
M}_{G,A}^{\bf d}({\bf C})$, satisfies $\rho(\delta)\,e^{2\pi i{\bf
d}/{\bf r}}=1$, where $\rho$ is the  natural homomorphism of $\pi_1(G/A)$
onto $A\i ({\bf G}_m)^s$ and  $e^{2\pi i{\bf d}/{\bf r}}$ stands for the
element $(e^{2\pi id_1/r_1},\ldots ,e^{2\pi id_s/r_s})$ of
$ ({\bf G}_m)^s$.
We denote by ${\cal
M}_{G,A}^\delta$ the open and closed substack $\pi^{-1}({\cal
M}_{G/A}^\delta)$ of ${\cal M}_{G,A}^{\bf d}$, where ${\bf
d}=(d_1,\ldots,d_s)$ is the unique element of ${\bf Z}^s$ such that $0\leq
d_t<r_t$ and  $\rho(\delta)\,e^{2\pi i{\bf d}/{\bf r}}=1$ (if $G$
is simply connected, $\rho$
 is bijective and
${\cal M}_{G,A}^\delta$ is simply ${\cal M}_{G,A}^{\bf d}$). The induced
morphism $\pi:{\cal M}_{G,A}^\delta \rightarrow {\cal M}_{G/A}^\delta$ is
surjective.

\ind  We will be mostly interested in the case when $A$ is the
center of $G$; then we will denote simply by  ${\cal M}_G^\delta$ the stack
${\cal M}^\delta_{G,A}$, for any choice of the isomorphism $ A\iso
\pprod_{j=1}^s\Mu_{r_j}$  (up to a canonical isomorphism, the stack ${\cal
M}^\delta_{G,A}$  does not depend on this choice).
If $\delta$ belongs to
$\pi_1(G)\subset\pi_1(G_{\rm ad})$, one gets $\rho(\delta)=1$ hence ${\bf
d}=0$: by the above remark, the notation ${\cal M}_{G}^\delta$ is thus
coherent with the one introduced in Prop.\ \ref{Bun_G}.

\medskip
\rem{Examples}\label{ex-M_G^\delta}
{\it a}) We take $G=SL_r$, $A=\Mu_r$.
 The group $C_AG$ is canonically isomorphic to
${\bf GL}_r$; the stack ${\cal M}_{{\bf SL}_r}^d$ can be identified with the
stack of
vector bundles $E$ on $X$ with an isomorphism
$\wedge^r E\iso {\cal O}_X(dp)$.

\ind {\it b}) We take for $G$ the group ${\bf O}_{2l}$ or ${\bf Sp}_{2l}$,
for $A$ its center, with the unique isomorphism $A\iso\Mu_2$. The group
$C_AG$ is the group $C{\bf O}_{2l}$ or $C{\bf Sp}_{2l}$ of automorphisms of
${\bf C}^{2l}$  respecting the bilinear form up to a (fixed) scalar. The
stack ${\cal M}_G^d$  can therefore be viewed as parameterizing vector
bundles $E$ on $X$ with a (symmetric or alternate) non-degenerate bilinear
form  with values in ${\cal O}_X(dp)$. Similarly, the stack
${\cal M}_{{\bf SO}_{2l}}^d$ parameterizes vector bundles $E$ on $X$ with a
non-degenerate quadratic  form $q:\sym^2E\rightarrow {\cal O}_X(dp)$ and an
{\it
orientation}, i.e.\ an isomorphism $\omega :\det E\iso {\cal O}_X(dlp)$ such
that
$\omega ^{\otimes2}$ coincides with the quadratic form induced by $q$ on
$\det E$.

\ind {\it c}) We take $G={\bf Spin}_r$, $A=\Mu_2$. Then $C_AG$ is the Clifford
group and ${\cal M}_{G,A}^{-1}$ is the moduli stack ${\cal M}^-_{{{\bf
Spin}_r}}$
considered in [O].

\medskip
\subsection\label{unif}  Choose
any element  $\zeta \in (LG_{\rm ad})^\delta({\bf C})$; reasoning as in Prop.\
\ref{Bun_G},  one gets a canonical isomorphism
${\cal M}_{G}^\delta\iso(\zeta^{-1}\,L_X{G}\,\zeta)\bk{\cal Q}_{\widetilde{G}}$
(see also [B-L], 3.6 for the case ${G}={\bf SL}_r$). In particular,
the stack ${\cal M}_{G}^\delta$ is connected.
Moreover, we see as in the proof of Prop.\ \ref{Pic-Bun} that the torsion
subgroup of
$\Pic({\cal M}_{G}^\delta)$ is canonically isomorphic to
$H^1(X,\pi_1(G)\,\widehat{}\ )$.

\ind Let us apply the above construction to the group
$\widetilde{G}$, with $A=\pi _1(G)$. Let $\delta\in \pi_1(G)$.
{}From the exact sequence (\ref{comp} {\it a}), we see that
$\zeta $ is the image of an element of $(L\widetilde{G})^\delta$. Comparing
with
Prop.\ \ref{Bun_G}, we see that the morphism $q_{G}^\delta:{\cal
Q}_{\widetilde{G}}\rightarrow {\cal M}_{G}^\delta$  factors as
$$q_{G}^\delta:{\cal Q}_{\widetilde{G}}\qfl{q_{\widetilde{G}}^\delta}
{\cal M}_{\widetilde{G}}^\delta \qfl{\pi} {\cal M}_{G}^\delta\ .$$
This shows us the way to determine the group $\Pic({\cal M}_{G}^\delta)$:
we will first compute
$\Pic({\cal M}_{G}^\delta)$ when $G$ is simply connected or $G={\bf SO}_{2l}$,
then
determine which powers of the generator(s) descend to ${\cal M}_{G}^\delta$.
\vskip1cm

\section{The Picard group of ${\cal M}_{{\bf PGL}_r}$}
 \ind
 According to (\ref{Bun_G}), the connected
components of ${\cal M}_{{\bf PGL}_r}$ are indexed by the integers $d$ with
$0\leq d<r$;  the component
${\cal M}_{{\bf PGL}_r}^d$ is dominated by the moduli stack ${\cal M}_{{\bf
SL}_r}^d$
parameterizing vector bundles $E$ on $X$ with an isomorphism $\wedge^rE\iso
{\cal
O}_X(dp)$ (\ref{ex-M_G^\delta} {\it a}).

\ind  Recall that  the {\it determinant
bundle} ${\cal D}$ on ${\cal M}_{{\bf SL}_r}^d$ is the line  bundle\break
$\det
R(pr_2)_*({\cal E})$, where ${\cal E}$ is the universal bundle on $X\times
{\cal M}_{{\bf SL}_r}^d$. It follows from  [B-L],
Prop.\ 9.2, that ${\cal D}$ generates $\Pic(
{\cal M}_{{\bf SL}_r}^d)$ and that its inverse
image on ${\cal Q}$ generates $\Pic({\cal Q})$. Therefore our problem is to
determine
which powers of  ${\cal D}$ descend to ${\cal
M}^d_{{\bf PGL}_r}$.

\th Proposition
\enonce The smallest power of ${\cal D}$ which descends to ${\cal M}_{{\bf
PGL}_r}^d$
is ${\cal D}^r$.
\endth\label{PGL_r}
{\it Proof}: Since it preserves the Killing form, the adjoint representation
defines a homomorphism ${\rm Ad}:{\bf GL}_r\rightarrow {\bf SO}_{r^2}$. Let
$f:{\cal
M}_{{\bf SL}_r}^d\rightarrow {\cal M}_{{\bf SO}_{r^2}}$ be the induced
morphism of stacks; since ${\rm Ad}$ factors through ${\bf PGL}_r$, $f$ factors
through ${\cal M}_{{\bf PGL}_r}^d$. By [L-S], the determinant bundle ${\cal
D}_{\bf
SO}$ on  ${\cal M}_{{\bf SO}_{r^2}}$ admits a square root ${\cal P}$;
one has  $f^*{\cal D}_{\bf SO}\cong{\cal D}^{2r}$ since the Dynkin index of
${\rm
Ad}$ is $2r$, hence $f^*{\cal P}\cong {\cal D}^r$, which implies that ${\cal
D}^r$
descends.

\ind  Let $J$ be the Jacobian of $X$, and
${\cal L}$  the Poincar\'e bundle on $X\times J$ whose restriction to
$\{p\}\times J$ is trivial. Consider the vector bundles
$${\cal F}={\cal L}^{\oplus (r-1)}\oplus {\cal L}^{1-r}(dp)\qquad {\rm
and}\qquad
{\cal G}= {\cal O}_X^{\oplus (r-1)}\oplus {\cal L}^{-1}(dp)$$ on $X\times J$.
We
denote by  $r^{}_J$ the multiplication by $r$ in $J$, and put $r^{}_{X\times
J}=\Id^{}_X\times  r^{}_J$. Since $r_{X\times J}^*\,{\cal L}\cong{\cal L}^r$,
one has
$r_{X\times J}^*\,{\cal G}\cong {\cal F}\otimes{\cal L}^{-1}$, hence the
projective
bundles $P({\cal F})$ and $r_{X\times J}^*\,P({\cal G})$ are isomorphic.
Therefore we
have a commutative\note{1}{By this we always mean 2-commutative, e.g.\ in our
case
the two functors $\pi\rond f$ and $g\rond r_J^{}$ are
isomorphic.}diagram of stacks
$$\diagramme{J & \phfl{f}{} & {\cal M}^d_{{\bf SL}_r}
&\cr
      \vfl{r^{}_J}{} &            &  \vfl{}{\pi}&\cr
             J &\phfl{g}{} &{\cal M}_{{\bf PGL}_r}^d &\kern-12pt ,
}$$
where $f$ and $g$ are the morphisms associated to ${\cal F}$ and $P({\cal G})$
respectively.

\ind Thus if ${\cal D}^k$
descends to ${\cal M}_{{\bf PGL}_r}^d$, the class of $f^*{\cal D}^k$ in the
N\'eron-Severi
group $NS(J)$ must be divisible by $r^2$. An easy computation shows that the
class of
$f^*{\cal D}$ in $NS(J)$ is $r(r-1)$ times the principal polarization; it
follows
that $r^2$ must divide $kr(r-1)$, which means that $r$ must divide $k$.\cqfd

\bigskip
\rem{Remark}  One can consider more generally the group $G={\bf SL}_r/\Mu_s$,
for each
integer $s$ dividing $r$, and the corresponding stacks ${\cal M}_G^d$ for $d\in
{r\over s}{\bf Z}$ (mod.$\,r{\bf Z}$). We can prove that {\it the line bundle
${\cal
D}^k$ descends to ${\cal M}_G^d$ if and only if $k$ is a multiple of
$s/(s,{r\over
s})$}. The ``only if" part is proved exactly as above, but the other
implication requires some descent theory on stacks which  lies beyond
the scope of this paper.

\vskip1cm
\section{The Picard group of ${\cal M}_{{\bf
PSp}_{2l}}$}
\ind According to Prop.\ \ref{Bun_G} the moduli stack ${\cal M}_{{\bf
PSp}_{2l}}$ has 2 components ${\cal M}_{{\bf PSp}_{2l}}^d$ $(d=0,1)$; the
component
${\cal M}_{{\bf PSp}_{2l}}^d$ is dominated by the algebraic stack ${\cal
M}_{{\bf
Sp}_{2l}}^d$ parameterizing vector bundles of rank $2l$ on $X$ with a
symplectic form
$\wedge^2E\rightarrow {\cal O}_X(dp)$ (Example \ref{ex-M_G^\delta} {\it b}).
Let
${\cal D}$ denote the determinant bundle on ${\cal M}_{{\bf Sp}_{2l}}^d$ (i.e.\
the
determinant of the cohomology of  the universal bundle on $X\times {\cal
M}_{{\bf
Sp}_{2l}}^d$); it is  the pull back of the determinant bundle ${\cal D}_0$ on
${\cal
M}_{{\bf SL}_{2l}}^d$ by the  $f:{\cal M}_{{\bf
Sp}_{2l}}^d\rightarrow {\cal M}_{{\bf SL}_{2l}}^d$ associated to the standard
representation.

  \th Lemma
\enonce The
group  $\Pic({\cal M}_{{\bf Sp}_{2l}}^d)$ is generated by ${\cal D}$.
\endth
\label{Sp^d}
{\it Proof}:  Consider the commutative diagram
$$\diagramme{
{\cal Q}_{{\bf Sp}_{2l}}&\hfl{F}{}& {\cal Q}_{{\bf SL}_{2l}}&\cr
\vfl{q_{{\bf Sp}_{2l}}^d}{} & & \vfl{}{q_{{\bf SL}_{2l}}^d}&\cr
{\cal M}_{{\bf Sp}_{2l}}^d & \hfl{f}{} & {\cal M}_{{\bf SL}_{2l}}^d& \kern-12pt
,
}$$
where $f$ and $F$ are induced by the embedding ${\bf Sp}_{2l}\rightarrow {\bf
SL}_{2l}$, and $q^d_G:{\cal Q}_G\rightarrow {\cal M}_G^d$ is the canonical
projection
(\ref{unif}). One has
 ${\cal D}=f^*{\cal D}_0$, $(q_{{\bf SL}_{2l}}^d)^*{\cal D}_0={\cal
O}_{{\cal Q}_{ {\bf SL}_{2l}}}(1)$ by [B-L],  5.5, and $F^*{\cal O}_{{\cal Q}_{
{\bf SL}_{2l}}}(1)={\cal O}_{{\cal Q}_{{\bf Sp}_{2l}}}(1)$ since the Dynkin
index of the standard representation of ${\bf Sp}_{2l}$ is $1$ ([L-S], Lemma
6.8). It follows that the homomorphism $(q_{{\bf Sp}_{2l}}^d)^*:\Pic({\cal M}_{
{\bf Sp}_{2l}}^d)\rightarrow \Pic({\cal Q}_{{\bf Sp}_{2l}})={\bf Z}\,{\cal
O}_{\cal
Q}(1)$ is surjective. On the other hand, the proof of Prop.\ 6.2
in [L-S] shows that it is injective; our assertion follows.\cqfd

\medskip
\ind In view of the above remarks, Prop.\ \ref{Pic-Bun} and (\ref{unif})
provide us
with an exact sequence $$0\rightarrow J_2\rightarrow \Pic({\cal M}_{{\bf
PSp}_{2l}}^d)\qfl{\pi^*}  \Pic({\cal M}_{{\bf Sp}_{2l}}^d)={\bf Z}\,{\cal D}\
;$$
 we now determine the image of $\pi^*$:
\th Proposition
\enonce  The smallest power of ${\cal D}$ which descends to ${\cal M}_{{\bf
PSp}_{2l}}^d$ is ${\cal D}$ if $l$ is even, ${\cal D}^2$ if $l$ is odd.
\endth \label{PSp}
{\it Proof}: The stack ${\cal M}_{{\bf Sp}_{2l}}^d$ parameterizes vector
bundles $E$
with a symplectic form $\varphi :\wedge^2E\rightarrow {\cal
O}_X(dp)$ (\ref{ex-M_G^\delta} {\it b}).  For such a pair, the form
$\wedge^2\varphi $ defines a quadratic form on $\wedge^2E$ with values in
${\cal
O}_X(2dp)$, hence an ${\cal O}_X$\tx valued quadratic form on $\wedge^2E(-dp)$.
Put
$N=l(2l-1)$; let $f_d:{\cal M}_{{\bf Sp}_{2l}}^d\rightarrow {\cal M}_{{\bf
SO}_{N}}$ be the morphism of stacks which associates to $(E,\varphi )$ the pair
$(\wedge^2E(-dp),\wedge^2\varphi )$. Since the representation $\wedge^2:{\bf
Sp}_{2l}\rightarrow {\bf SO}_N$ factors through ${\bf PSp}_{2l}$, the morphism
$f_d$ factors as
$$f_d: {\cal M}_{{\bf Sp}_{2l}}^d\rightarrow {\cal M}_{{\bf
PSp}_{2l}}^d\rightarrow {\cal M}_{{\bf SO}_{N}}\ .$$  The pull back under $f_d$
of
the determinant bundle on ${\cal M}_{{\bf SO}_N}$ is ${\cal D}^{2l-2}$ ($2l-2$
is the
Dynkin index of the representation $\wedge^2$). But we know by [L-S] that this
determinant bundle admits a square root, hence ${\cal D}^{l-1}$
descends to ${\cal M}_{{\bf
PSp}_{2l}}^d$. On the other hand, the same argument applied to the adjoint
representation shows that ${\cal D}^{2l}$ descends (see the proof of Prop.\
\ref{PGL_r}). We conclude that ${\cal D}^2$ descends, and that ${\cal D}$
descends
when $l$ is even.

\ind To
prove that ${\cal D}$ does not descend when $l$ is odd, we use the notation of
the proof of Prop.\
\ref{PGL_r}, and consider on $X\times J$ the vector bundle ${\cal H}={\cal
L}^{\oplus
l}\oplus {\cal L}^{-1}(dp)^{\oplus l}$, endowed with the standard hyperbolic
alternate form with values in ${\cal O}(dp)$. We see as in {\it loc.\ cit.}
that the ${\bf PSp}_{2l}$\tx bundle associated to ${\cal H}$ descends under the
isogeny $2^{}_J$ (observe that ${\cal H}\otimes{\cal L}$ descends, and use the
exact sequence $\ 1\rightarrow {\bf G}_m\rightarrow C{\bf Sp}_{2l}\rightarrow
$ $\rightarrow {\bf PSp}_{2l}\rightarrow 1$). Therefore the morphism
$h:J\rightarrow {\cal M}_{{\bf Sp}_{2l}}^d$ defined by ${\cal H}$ fits in a
commutative diagram  $$\diagramme{J &
\phfl{h}{} & {\cal M}^d_{{\bf Sp}_{2l}} &\cr
      \vfl{2_J}{} &            &  \vfl{}{}&\cr
             J &\phfl{}{} &{\cal M}^d_{{\bf PSp}_{2l}} &\kern-12pt .
}$$
Since the class of $f^*{\cal D}$ in $NS(J)$ is $2l$ times the principal
polarization, it follows that ${\cal D}$ does not descend.\cqfd

  \vskip1cm

\section{The Picard group of ${\cal M}_{{\bf
PSO}_{2l}}$}
 \subsection Let us consider first the moduli stack ${\cal M}_{{\bf SO}_r}$,
for $r\geq 3$.
It has two components ${\cal M}_{{\bf SO}_r}^w$, distinguished by the second
Stiefel-Whitney class $w\in \Mu_2$. The Picard group of these stacks
is essentially described in [L-S]: to each theta-characteristic  $\kappa $
on $X$ is associated a
Pfaffian line bundle ${\cal P}_\kappa $ whose square is the determinant
bundle ${\cal D}$ (determinant of the cohomology of
the universal bundle on $X\times {\cal M}_{{\bf SO}_{r}}^w$); according to
Prop.\
\ref{Pic-Bun}, there is a canonical exact sequence
$$0\rightarrow J_2\qfl{\lambda  } \Pic({\cal
M}_{{\bf SO}_r}^w) \longrightarrow {\bf Z}\rightarrow 0\ ,$$
where
 the torsion free quotient
is  generated by any of the   ${\cal P}_\kappa $'s.
\medskip
 \ind  We can actually be  more
precise. Let $\theta (X)$ be the subgroup of  $\Pic(X)$ generated by the
theta-characteristics; it is  an extension of ${\bf Z}$ by $J_2$.
\th Proposition
\enonce  The map $\kappa \mapsto{\cal P}_\kappa $ extends by
linearity to an isomorphism  ${\cal P}:\theta (X)\iso \Pic({\cal M}_{{\bf
SO}_r}^w)$,
which coincides with $\lambda $ on $J_2$.
\endth\label{\theta }
\ind In other words, we have a canonical isomorphism of extensions
$$\diagramme{
0\rightarrow & J_2 & \phfl{}{} & \theta (X) & \phfl{}{}&{\bf Z}
&\rightarrow 0 &\cr
&\left \| \vbox to 6truemm{}\right. &  & \vfl{\cal P}{} & & \left \| \vbox to
6truemm{}\right. &&\cr
0\rightarrow & J_2 & \phfl{\lambda }{} & \Pic({\cal M}_{{\bf
SO}_{r}}^w)  & \phfl{}{}&{\bf Z} &\rightarrow 0 &\kern-12pt .
 }$$

{\it Proof}: It suffices to prove the formula ${\cal
P}_{\kappa\otimes\alpha }={\cal P}_{\kappa}\otimes \lambda (\alpha)$ for any
theta-characteristic $\kappa $ and element $\alpha $ of $J_2$.

\ind Let ${\cal L}$ be the Poincar\'e bundle on $X\times J$, normalized so that
 its
restriction to $\{p\}\times J$ is trivial. Put $d=0$ if $w=1$, $d=1$ if $w=-1$.
The
vector bundle  ${\cal
L}(dp)\oplus {\cal L}^{-1}(-dp)\oplus  {\cal O}^{r-2}$, with its natural
quadratic
form and orientation, defines a morphism $g:J\rightarrow {\cal M}_{{\bf
SO}_r}^w$.
Let us identify $J$ with $\Pic^{\rm o}(J)$   via the
principal polarization. Then the required formula is a consequence of  the
following
two assertions:  \indp {\it a}) One has $g^*{\cal
P}_{\kappa\otimes\alpha }=(g^*{\cal P}_{\kappa})\otimes \alpha $ for every
theta-characteristic $\kappa $ and element $\alpha $ of $J_2$; \indp{\it b})
The map
$g^*:\Pic({\cal M}_{{\bf SO}_r}^w)_{tors}\rightarrow J_2$ is the inverse
isomorphism of $\lambda  $. \ind Let us prove {\it a}). The line bundle
$g^*{\cal
P}_\kappa $  is the pfaffian bundle associated to the quadratic bundle  ${\cal
L}(dp)
\oplus {\cal L}^{-1}(-dp)$) and to $\kappa $. Now it follows from the
construction in
[L-S] that for any vector bundle $E$ on $X\times S$, the pfaffian of the
cohomology of
$E\oplus (K_X\otimes E^*)$, endowed with the standard hyperbolic form with
values in
$K_X$, is the determinant of the cohomology of $E$. Because the choice of
${\cal L}$
ensures that the determinant of the cohomology is the same for ${\cal L}$ and
${\cal
L}(p)$, we conclude that $g^*{\cal
P}_\kappa $ is the determinant of the cohomology of ${\cal L}\otimes\kappa $,
i.e.\
the  line bundle  ${\cal O}_J(\Theta _\kappa )$. Since $\Theta
_{\kappa\otimes\alpha
}=\Theta _\kappa +\alpha  $, the assertion {\it a}) follows.
\ind Since we already know
 that $\Pic({\cal M}_{{\bf
SO}_r}^\pm)_{tors}$ is isomorphic to $J_2$ (Prop.~\ref{Pic-Bun}), {\it a})
implies that $g^*$ is surjective, and therefore  bijective. Hence $u=g^*\rond
\lambda $ is an automorphism of $J_2$. This construction extends to any family
of curves $f:{\cal X}\rightarrow {\cal S}$, defining an automorphism of the
local system $R^1f_*(\Mu_2)$ over ${\cal S}$.   Since the
mono\-dromy group of this local system  is the full symplectic
group ${\bf Sp}(J_2)$  for the universal family of
curves, it follows that $u$ is the identity.\cqfd

 \medskip
\subsection
This settles the case of the group ${\bf SO}_r$; let us now consider the group
${\bf PSO}_r$, for  $r=2l\geq 4$. The moduli space ${\cal M}_{{\bf PSO}_{2l}}$
has $4$
components, indexed by the center $Z$ of ${\bf Spin}_{2l}$. This group consists
of the
elements $\{1,-1,\varepsilon ,-\varepsilon \}$ of the Clifford algebra $C({\bf
C}^{2l})$,  with $\varepsilon ^2=(-1)^l$ ([Bo], Alg\`ebre IX). Each component
${\cal
M}_{{\bf PSO}_{2l}}^\delta$, for $\delta\in Z$, is dominated by the algebraic
stack
${\cal M}_{{\bf SO}_{2l}}^\delta$ (\ref{M_G^\delta}). For
$\delta\in\{\pm 1\}$, this is the same stack as above; the stack ${\cal
M}_{{\bf
SO}_{2l} }^\varepsilon  \cup {\cal M}_{{\bf SO}_{2l}}^{-\varepsilon }$
parameterizes vector bundles with a quadratic form with values in ${\cal
O}_X(p)$ and an orientation (\ref{ex-M_G^\delta} {\it b}). Changing the sign of
the orientation exchanges the two components  ${\cal M}^\varepsilon $ and
${\cal
M}^{-\varepsilon }$ (this corresponds to the fact that $\varepsilon $ and
$-\varepsilon $ are exchanged by the outer automorphism of ${\bf Spin}(2l)$
defined by conjugation by an odd degree element of the Clifford group).
\label{PSO}
\medskip
\th Lemma
\enonce The torsion free quotient of $\Pic({\cal M}_{{\bf
SO}_{2l}}^{\pm\varepsilon})$ is generated by the determinant bundle ${\cal D}$.
\endth
{\it Proof}: The same proof as in Lemma \ref{Sp^d} shows that the pull back of
${\cal
D}$ by the morphism $q_{{\bf SO}_{2l}}^{\pm\varepsilon }:{\cal Q}_{{\bf\rm
Spin}_{2l}}\rightarrow {\cal M}_{{\bf SO}_{2l}}^{\pm\varepsilon}$ is ${\cal
O}_{\cal Q}(2)$ (the Dynkin index of the standard representation of
${\bf SO}_{2l}$ is $2$). Therefore it suffices to prove that ${\cal D}$ has no
square root in $\Pic({\cal M}_{{\bf SO}_{2l}}^{\pm\varepsilon})$.
\ind   Let $V$ be a $l$\tx dimensional vector space; we consider the
vector bundle\break $T =(V\otimes_{\bf C}{\cal O}_X)\oplus (V^*\otimes_{\bf
C}{\cal O}_X(p))$, with the obvious hyperbolic quadratic form\break
$q:\sym^2T\rightarrow {\cal O}_X(p)$ and  isomorphism $\omega :\det T \iso
{\cal
O}_X(lp)$. We choose the sign of $\omega $ so that the triple
$T^\varepsilon:=(T,q,\omega ) $ defines a rational point of ${\cal M}_{{\bf
SO}_{2l} }^{\varepsilon}$, and put $T^{-\varepsilon }:=(T,q,-\omega )\in {\cal
M}_{{\bf SO}_{2l} }^{-\varepsilon}({\bf C})$. The group $G={\bf GL}(V)$ acts on
$T$, and this action  preserves the quadratic form and the orientation. This
defines a morphism $\iota$ of the stack $BG$ classifying $G$\tx torsors into
${\cal M}_{{\bf SO}_{2l}}^{\pm\varepsilon}$: if $S$ is a ${\bf C}$\tx scheme
and
$P$ a $G$\tx torsor on  $S$, one puts $\iota (P)=P\times^GT^{\pm\varepsilon
}_S$. \ind Recall [L-MB] that  the ${\bf C}$\tx stack $BG$  is the quotient of
$\Spec {\bf C}$ by the trivial action of $G$; in particular, line bundles on
$BG$ correspond in a one-to-one way to $G$\tx linearizations of the trivial
line
bundle on $\Spec {\bf C}$, that is to characters of $G$. In our situation, the
line bundle $\iota ^*{\cal D}$ will correspond to the character of $G$ by which
$G$ acts on $\det R\Gamma (X,T)$. As $G$\tx modules, we have
$$\det R\Gamma (X,T)\cong\det R\Gamma (X,V\otimes_{\bf C}
{\cal O}_X)\otimes \det R\Gamma (X,V^*\otimes_{\bf C}{\cal O}_X(p))\ .$$ Now if
$L$ is a line bundle on $X$, the $G$\tx module $\det R\Gamma (X,V\otimes_{\bf
C}L)$ is  isomorphic to $\det(V\otimes H^0(L))\otimes \det(V\otimes
H^1(L))^{-1}=\det(V)^{\chi (L)}$. We conclude that $\det R\Gamma (X,T)$ is
isomorphic to $\det(V^*)$, i.e.\ that $\iota ^*{\cal D}$ corresponds to the
character $\det^{-1}:G\rightarrow {\bf C}^*$. Since $\det$ generates
$\Hom(G,{\bf C}^*)$, our assertion follows.\note{1}{This argument has been
shown
to us by V. Drinfeld.}\cqfd

\th Proposition
\enonce Let $\delta\in Z$.
The line bundle ${\cal D}$ {\rm (}resp.\ ${\cal D}^2)$ descends on ${\cal
M}_{{\bf PSO}_{2l}}^\delta$ if $l$ is even {\rm (}resp.\ odd{\rm );} the
corresponding line bundles on ${\cal M}_{{\bf PSO}_{2l}}^\delta$ generate the
 Picard group. \endth
{\it Proof}: We first prove that the Pfaffian bundles ${\cal P}_\kappa $ do not
descend to ${\cal M}_{{\bf PSO}_{2l}}^\delta$. If $\delta\in\{\pm\varepsilon
\}$, this
follows from the above lemma. If $\delta\in\{\pm1\}$, we consider the action of
$J_2$
on ${\cal M}_{{\bf SO}_{2l}}^\delta$ deduced from the embedding $\Mu_2\i {\bf
SO}_{2l}$: each element $\alpha \in J_2$ (trivialized at $p$) defines an
automorphism
-- still denoted $\alpha $ -- of the stack ${\cal M}_{{\bf SO}_{2l}}^\pm$,
which
maps a quadratic bundle $(E,q,\omega )$ onto $(E\otimes\alpha, q\otimes
i_\alpha
,\omega \otimes i_\alpha ^{\otimes l}) $, where $i_\alpha :\alpha ^2\iso {\cal
O}_X$
is the isomorphism which coincides at $p$ with the square of the given
trivialization.
\ind We claim that $\alpha ^*{\cal P}_\kappa $ {\it is isomorphic to} ${\cal
P}_{\kappa
\otimes\alpha }$ for every theta-characteristic $\kappa $ and element $\alpha$
of
$J_2$. This is easily seen by using the following characterization of  ${\cal
P}_\kappa $ ([L-S], 7.10): let ${\cal E}$ be the universal bundle on $X\times
{\cal
M}_{{\bf SO}_{2l}}^\pm$; then the divisor $\Theta _\kappa :=\div Rpr_{2*}({\cal
E}\otimes\kappa )$ is divisible by $2$ in  $\Div{\cal M}_{{\bf SO}_{2l}}^\pm$,
and
${\cal P}_\kappa $ is the line bundle associated to ${1\over 2}\Theta _\kappa
$.
By construction $(1_X\times
\alpha )^*{\cal E}$ is isomorphic to ${\cal E}\otimes\alpha $, hence
$$\alpha ^*\Theta _\kappa =\div Rpr_{2*}((1_X\times
\alpha )^*{\cal E}\otimes\kappa )=\div Rpr_{2*}({\cal E}\otimes\alpha
\otimes\kappa
)=\Theta _{\kappa \otimes\alpha }\ ,$$which implies our claim. Since the map
$\kappa
\mapsto {\cal P}_\kappa $ is injective (Prop.\ \ref{\theta }), we conclude that
${\cal
P}_\kappa $ does not descend.
\ind The rest of the proof follows closely the symplectic case (Prop.\
\ref{PSp}).
For $d=0,1$, the representation $\wedge^2$ defines a morphism of stacks
 $g_d:{\cal M}_{{\bf SO}_{2l}}^d\rightarrow {\cal M}_{{\bf
SO}_{N}}$, which factors through
${\cal M}_{{\bf PSO}_{2l}}^d$.  The pull back under $g_d$ of a square root of
the determinant bundle is ${\cal D}^{l-1}$; since
 ${\cal D}^{2l}$ descends, one concludes
 that ${\cal D}$ descends when $l$ is even and ${\cal D}^2$ when $l$ is odd.
\ind To prove that ${\cal D}$ does not descend when $l$ is odd, one considers
the
quadratic bundle ${\cal H}^\delta$ on $X\times J$ defined by
$$\eqalign{
{\cal H}^\delta&={\cal L}^{\oplus l}\oplus({\cal L}^{-1})^{\oplus l}\quad
\hbox{if }
\delta=1 \cr
&= {\cal L}(p)^{\oplus l}\oplus{\cal L}^{-1}(-p))^{\oplus l}\quad \hbox{if }
\delta=-1 \cr
&= ({\cal L}\oplus{\cal L}^{-1}(p))^{\oplus l}\quad \hbox{if }
\delta=\pm\varepsilon \ ,}$$
with the standard hyperbolic quadratic form, and opposite orientations for the
cases
$\delta=\varepsilon $ and $\delta=-\varepsilon $.
\ind As above, this gives rise to a commutative diagram
$$\diagramme{J &
\phfl{h}{} & {\cal M}^\delta_{{\bf SO}_{2l}} &\cr
      \vfl{2_J}{} &            &  \vfl{}{}&\cr
             J &\phfl{}{} &{\cal M}^\delta_{{\bf PSO}_{2l}} &\kern-12pt ,
}$$from which one deduces that ${\cal D}$ does not descend, since the class of
$h^*{\cal D}$ in $NS(J)$ is $2l$ times the principal polarization.\cqfd

\vskip2cm
\centerline{\gros Part II: The Picard group of the moduli space}
\vskip1cm

\section{${\bf C}^*$\tx extension associated to group actions}\label{ext}

\ind  This part is devoted to the Picard group of the
moduli space $M_G$. The case of a simply connected group being known, we will
consider
$M_G$ as a quotient of  $M_{\widetilde{G}}$ by the finite group
$H^1(X,\pi_1(G))$.
Therefore we will first develop some general tools to study the Picard group of
a
quotient variety.
\subsection Let $H$ be a finite group
acting on an integral projective variety $M$ over ${\bf C}$ (or, more
generally, any variety with $H^0(M,{\cal O}_M^*)={\bf C}^*$), and $L$  a
line bundle on $M$ such that $h^*L\cong L$ for all $h\in H$. We associate to
this situation a central ${\bf C}^*$\tx  extension  $$1\rightarrow {\bf
C}^*\longrightarrow {\cal E}(H,L)
 \qfl{p} H\rightarrow 1 \ ;$$
the group ${\cal E}(H,L)$ consists of pairs $(h,\tilde h)$, where $h$ runs over
$H$
and $\tilde h$ is an automorphism of $L$ covering $h$, and $p$ is the first
projection.

\smallskip

\subsection We will need a few elementary properties of these groups:
\label{list}
\ind {\it a}) Let $f:M'\rightarrow M$ be a $H$\tx equivariant rational map.
Pulling back automorphisms induces an isomorphism $f^*:{\cal
E}(H,L)\rightarrow {\cal E}(H,f^*L)$.

\smallskip
\ind  {\it b}) Recall that the isomorphism classes of central ${\bf C}^*$\tx
extensions of $H$ form a group, canonically isomorphic
to $H^2(H,{\bf C}^*)$.
Let $r$ be a positive integer. The extension ${\cal E}(H,L^r)$ is isomorphic
to the sum of $r$ copies of ${\cal E}(H,L)$; more precisely, the homomorphism
$\varphi_r :{\cal E}(H,L)\rightarrow {\cal E}(H,L^r)$ given by $\varphi_r
(h,\tilde h)= (h,\tilde h^{\otimes r})$ is a surjective homomorphism, with
kernel
the group $\Mu_r$ of $r$\tx roots of unity, and therefore induces an
isomorphism
of  the push forward of ${\cal E}(H,L)$ by the endomorphism $t\mapsto t^r$ of
${\bf C}^*$ onto  ${\cal E}(H,L^r)$.
\smallskip
\ind {\it c}) Let $M'$ be another projective variety, $H'$ a finite group
acting
on $M'$, $L'$ a line bundle on $M'$ preserved by $H'$. The map $ {\cal E}(H
,L)\times {\cal E}(H' ,L') \rightarrow {\cal E}(H\times H'
,$ $L\boxtimes L')$ given by $((h,\tilde h),(h',\tilde h'))\mapsto (h\times
h',\tilde h\,{\scriptstyle\boxtimes}\,\tilde h')$ is a surjective
homomorphism, with kernel ${\bf C}^*$ embedded by $t\mapsto(t,t^{-1} )$.
\smallskip
\ind {\it d})  Let $K$ be a normal subgroup of $H$. The group  $H/K$ acts
on $M/K$; let $L_0$ be a line bundle on $M/K$ preserved by this action, and $L$
the pull back of $L_0$ to $L$. Then the extension ${\cal E}(H ,L)$ is the
 pull back of ${\cal E}(H/K ,L_0)$ by the projection $H\rightarrow H/K
$.\medskip

\subsection A $H$\tx {\it linearization} of $L$ is a section of the
extension ${\cal E}(H,L)$.  Such a linearization allows us to define,
for each point $m$ of $M$, an action of
the stabilizer $H_m$ of $m$ in $H$ on the fibre  $L_m$;  this action is given
by a
character of $H_m$, denoted by $\chi _m$.
\label{char}
\ind Let $\pi:M\rightarrow M/H$ be the quotient map; if $L_0$ is a line bundle
on
the quotient $M/H$, the line bundle $L=\pi^*L_0$ has a canonical $H$\tx
linearization.  By construction it has the property that at each point $m$ of
$M$, the character $\chi _m$ of  $H_m$  is {\it trivial}.
The converse is true (``Kempf's lemma"),  and is actually quite easy to prove
in
our set-up. Since two $H$\tx linearizations differ by a character of $H$, we
can
restate this lemma as follows: assume that  $L$ admits a  $H$\tx
linearization; then $L$ {\it descends to
$M/H$ if and
only if there exists a character $\chi $ of $H$ such that $\chi _m=\chi
_{|H_m}$
for all} $m\in M$.
\ind It follows from the above description that the kernel of the homomorphism
$\pi^*: \Pic(M/H)\rightarrow \Pic(M)$ consists of the $H$\tx linearizations of
${\cal O}_M$ such that the associated characters $\chi _m$ are trivial, i.e.\
of
the characters of $H$ which are trivial on the stabilizers $H_m$. In
particular,
if the subgroups $H_m$ generate $H$, $\pi^*$ is injective.
 \subsection \label{produit} Let $M'$ be another projective variety with an
action of $H$, and $L'$ a line bundle  admitting a  $H$\tx
linearization. The   $H$\tx linearizations of $L$ and $L'$ define a  $H$\tx
linearization of $L\boxtimes L'$; at each point $(m,m')$, the corresponding
character of $H_{(m,m')}=H_m\cap H_{m'}$ is the product of the characters $\chi
_m$ of $H_m$ and $\chi' _{m'}$ of $H_{m'}$ associated to the linearizations of
$L$ and $L'$. As a consequence, assume that $L$ and $L\boxtimes L'$ descend to
$M/H$ and $(M\times M')/H$ respectively, and  that $H=\cup_m H_m$; {\it then
$L'$
descends to} $M'/H$.
\medskip
\subsection\label{section-mod}
 From (\ref{list} {\it b})
  we see that {\it  the smallest positive integer $n$ such that $L^n$ admits a
$H$\tx linearization is the order of  ${\cal E}(H,L)$ in $H^2(H,{\bf
C}^*)$}. We want to know which powers of $L^n$ descend to $M/H$.

 \ind Let $r$ be a multiple of $n$. The class  $e$ of ${\cal E}(H,L)$ in
$H^2(H,{\bf C}^*)$ comes from an element of $H^2(H,\Mu_r)$, which means that
there
exists a cocycle $c\in Z^2(H,\Mu_r)$ representing $e$, or in other words a map
$\sigma :H\rightarrow {\cal E}(H,L)$ such that $p\rond\sigma ={\rm Id}_H$ and
$\sigma
(hh')\equiv \sigma (h)\,\sigma (h')$ (mod.$\,\Mu_r)$ -- let us call such a map
a
{\it section}  ({\it
 mod.}~$\Mu_r)$ {\it of} ${\cal E}(H,L)$. Composing $\sigma $ with the
homomorphism $\varphi _r: {\cal E}(H,L)\rightarrow {\cal E}(H,L^r)$ (\ref{list}
{\it b}) gives
 a section of the extension ${\cal E}(H,L^r)$, that is  a $H$\tx
linearization of $L^r$.

 \ind Let $m$ be a point of $M$. Using this $H$\tx linearization we
get a character $\chi _m$  of  $H_m$ (\ref{char}), which can be computed as
follows: for $h\in H_m$ the element $\sigma (h)$ acts on $L_m$,   and we have
$\chi _m(h)=(\sigma (h)_m)^r$. Assume moreover that $h^r=1$ for all $h\in H$;
then the element $\sigma (h)^r$ of ${\cal
E}(H,L)$ belongs to the center ${\bf C}^*$. Thus $L^r$ {\it endowed with the
$H$\tx linearization deduced from $\sigma $ descends to $M/H$ if and only if
$\sigma
(h)^r=1$ for all $h$ in} $\cup H_m$. Using \ref{char} we can conclude:

\th
 Proposition
\enonce  Assume that the order of ${\cal
E}(H,L)$ in $H^2(H,{\bf C}^*)$ and of every element of $H$ divides $r$.
 Let $\sigma :H\rightarrow {\cal E}(H,L)$ be a section {\rm (mod.}\
 $\Mu_r)$. Then $L^r$ descends to $M/H$ if and only if there exists a character
$\chi $ of $H$ such that $\sigma (h)^r=\chi (h)$ for all $h\in H$ fixing some
point of $M$.
\endth
\label{desc-prelim}
\smallskip
\ind  In the applications we have in mind we will always have $\cup H_m=H$. In
this case we get the following condition, which depends  only on the extension
${\cal E}(H,L)$ and not on the variety $M$:

\th Corollary
\enonce Assume that every element of $H$ fixes some point in
$M$. Then $L^r$ descends to $M/H$ if  and only if the map $h\mapsto \sigma
^r(h)$ from $H$ to ${\bf C}^*$ is a homomorphism.
\label{desc-cor}
\endth
\smallskip
\subsection From now on we will assume that the finite group $H$ is {\it
abelian}. In that case there is a canonical isomorphism of $H^2(H,{\bf C}^*)$
onto the group ${\rm Alt}(H,{\bf C}^*)$ of
bilinear alternate forms on $H$ with values in ${\bf C}^*$ (see e.g.\ [Br],
V.6, exer.\ 5) : to a central
${\bf C}^*$\tx
 extension $\widetilde{H}\qfl{p} H$ corresponds the form $e$ such that $\
e(p(x)
,p(y ))=$ $xyx^{-1}y ^{-1}\in\Ker p={\bf C}^*$. Conver\-sely,
given $e\in{\rm
Alt}(H,{\bf C}^*)$, one defines an
 extension of $H$ in the following way:  choose any bilinear form $\varphi
:H\times H\rightarrow {\bf C}^*$ such that $e(\alpha ,\beta )=\varphi (\alpha
,\beta )\varphi (\beta ,\alpha )^{-1}$; take $\widetilde{H}=H\times{\bf C}^*$,
with
the multiplication law given by $$(\alpha ,t)\,(\beta ,u)=(\alpha +\beta
,tu\,\varphi
(\alpha ,\beta ))\ ,$$  the homomorphism $p:\widetilde{H}\rightarrow H$ being
given by
the first projection.

\medskip
\subsection Let $r$ be an integer such that $rH=0$. Then the group $H^2(H,{\bf
C}^*)\cong$ ${\rm Alt(H,{\bf C}^*)}$ is
also annihilated by $r$. Let $e\in {\rm Alt}(H,{\bf C}^*)$; we can choose the
bilinear form $\varphi $ with values in $\Mu_r$. Consider the extension defined
as
above by  $\varphi $. The map $\sigma :H\rightarrow \widetilde{H}$ defined by
$\sigma (\alpha )=(\alpha ,1)$ is a section  (mod.\ $\Mu_r)$. An easy
computation
gives $\sigma (\alpha )^r= \varphi (\alpha ,\alpha )^{{1\over
2}r(r-1)}\in\{1,-1\}$. One has $\sigma ^{2r}(\alpha )=1$, and   $\sigma (\alpha
)^r=1$ for all $\alpha $ if $r$ is odd.  If $r$ is even, the function
$\varepsilon :\alpha \mapsto \sigma (\alpha )^r$ is ``quadratic" in the sense
that $\varepsilon (\alpha +\beta )=\varepsilon (\alpha )\varepsilon (\beta
)\,e(\alpha ,\beta )^{r/2}$. In particular, we see that $\sigma ^r$ is a
homomorphism if and only if the alternate form $e^{r/2}$ (with values in
$\Mu_2$)
is trivial.  In summary: \smallskip
\th Proposition
\enonce Assume $H$ is commutative, annihilated by $r$;
let $e$ be the alternate form associated to ${\cal E}(H,L)$. The line bundle
$L^{2r}$ descends to $M/H$; moreover $L^r$ descends,  except if  $r$ is even
and
the form $e^{r/2}$ is not trivial. In this last case, if every element of $H$
has
some fixed point on $M$, $L^r$ does not descend.
\label{desc-prop}
\endth
\smallskip
\rem{ Example}  \label{Theta } Let
$A$ be an abelian variety of dimension $g\ge 1$, and $\Theta $ a divisor on $A$
defining a principal polarization. Let $A_r$ be the kernel of the
multiplication by
$r$ in $A$. The  group ${\cal E}(A_r,{\cal O}(r\Theta ))$ is the Heisenberg
group
which plays a fundamental role in Mumford's theory of theta
functions; the corresponding alternate form $e_r:A_r\times A_r\rightarrow
\Mu_r$
is the Weil pairing. The group $A_r$ acts on the linear system $|r\Theta |$,
and the morphism $A\rightarrow |r\Theta |^*$ associated to this linear system
is
$A_r$\tx equivariant; therefore by (\ref{list} {\it a}), the extension  ${\cal
E}(A_r,{\cal O}_{|r\Theta |^*}(1))$ is  isomorphic to ${\cal
E}(A_r,{\cal O}(r\Theta ))$. It follows easily that ${\cal E}(A_r,{\cal
O}_{|r\Theta |}(1))$ corresponds to the nondegenerate form $e_r^{-1}$.
Let $s$ be a positive integer dividing $r$; an easy computation shows that\par

\subsec{\it  the restriction of $e_r$ to $A_s$ is equal to} $e_s^{r/s}$.
\label{restric} \ind We conclude from the proposition  that:\par  -- the line
bundle
${\cal O}(2s)$ descends to $|r\Theta |/A_s$; \par
-- the line bundle ${\cal O}(s)$ descends to $|r\Theta|/A_s$ if $s$ is odd or
$r/s$ is even, but does {\it not} descend if $s$ is even and $r/s$ odd.

\vskip1cm

\section{The moduli space $M_G$}
\subsection\label{coarse}
Recall [R1, R2] that  a $G$\tx
bundle $P$ on $X$ is {\it semi-stable} (resp.\ {\it stable}) if for every
parabolic subgroup $\Pi $,  every dominant  character
$\chi$ of $\Pi $ and  every $\Pi $\tx bundle $P'$ whose associated $G$\tx
bundle
is isomorphic to $P$, the line bundle $P'_\chi $ has nonpositive (resp.\
negative)
degree.
\ind  Let
$\rho:G\rightarrow  G^{\prime}$ be a homomorphism of semi-simple  groups, and
$P$
a  $G$\tx bundle. If $P$ is semi-stable  the $G'$\tx bundle $P_\rho=P\times^G
G^{\prime}$ is semi-stable;  the
converse is true if $\rho $ has finite kernel. In particular $P$ is
semi-stable if and only if its adjoint bundle ${\rm Ad}(P)$ is semi-stable.
\ind We denote by $M_G$ the coarse moduli space  of semi-stable
principal $G$\tx bundles on $X$ ({\it loc.\ cit.}). It is a projective
normal variety. Let ${\cal M}_G^{ss}$ be the open
substack of ${\cal M}_G$  corresponding to semi-stable $G$\tx bundles; there
is a canonical  surjective morphism  $f:{\cal M}_G^{ss}\rightarrow
M_G$.  For $\delta\in \pi_1(G)$, $f$ maps the component $({\cal
M}_G^{ss})^\delta$ onto the
connected component  $M_G^\delta$ of $M_G$ parameterizing $G$\tx bundles of
degree $\delta$.

\ind The definition of (semi-) stability extends to any reductive group $H$:
a
$H$\tx bundle $P$ is
semi-stable (resp.\ stable) if and only if the $(H/Z^{\rm o})$\tx bundle
$P/Z^{\rm
o}$ has the same property, where $Z^{\rm o}$ is the neutral
component of the center of $H$. The construction of the moduli space $M_H$
makes sense in this set-up; each component of $M_H$ is normal and projective.
\ind  Let $Z$ be the center of $G$; we choose  an isomorphism
$Z\iso \pprod_{}^{} \Mu_{r_j}$. Let $\delta\in\pi_1(G_{\rm ad})$.  The
construction of the ``twisted" moduli stack ${\cal M}_{G}^\delta$ (Section
\ref{twist}) obviously makes sense in the framework of coarse moduli
spaces. We get a coarse moduli space $M_{G}^\delta$, which
parameterizes semi-stable $C_ZG$\tx bundles with
determinant ${\cal O}_X({\bf d}p)$, such that the associated $G_{\rm ad}$\tx
bundle
has degree $\delta$, with $\rho (\delta)\,e^{2\pi i{\bf d}/{\bf r}}=1$
(\ref{deg}). The open substack ${\cal M}_{G}^{\delta, ss}$ of ${\cal
M}_{G}^\delta$
parameterizing semi-stable bundles maps surjectively onto
$M_{G}^\delta$.  If $A$ is a central subgroup of $G$, there is a canonical
morphism
$\pi:M_{G}^\delta\rightarrow  M_{G/A}^\delta$ which
is a (ramified) Galois covering with Galois group $H^1(X,A)$. The  next lemma
will allow us to compare the Picard groups of these moduli spaces by
 applying the results of section
\ref{ext}, in particular Prop.\ \ref{desc-prop}, to the action of $H^1(X,A)$
on $M_G^\delta$.

\th Lemma
\enonce Let $\delta\in \pi _1(G_{\rm ad})$.
\ind {\rm a)} The moduli space $M_{G}^\delta$ is unirational.
\ind {\rm b)} Any finite order automorphism  of $M_{G}^\delta$ has
fixed points.
\endth\label{unirat}
{\it Proof}: a) The proof in [K-N-R], Cor.\ 6.3, for the untwisted case
extends in a straightforward way: by (\ref{unif}) we have a surjective
morphism ${\cal Q}_{\widetilde{G}}\rightarrow {\cal M}_{G}^\delta$;
so the
open subset of ${\cal Q}_{\widetilde{G}}$ parameterizing semi-stable bundles
maps surjectively onto $M_G^\delta$.
Since ${\cal Q}_{\widetilde{G}}$ is a direct
limit
of an increasing sequence of generalized Schubert varieties, which are
rational, the lemma follows.

\ind b)  This is actually true for any finite order automorphism
$g$ of a  projective unirational variety $M$. One (rather sophisticated)
proof goes as follows: there exists a desingularization $\widetilde{M}$ of
$M$ to which $g$ lifts to an automorphism $\tilde g$ [H], necessarily
of finite order. Since $H^i(\widetilde{M},{\cal O}_{\widetilde{M}})$ is
zero for $i>0$, we deduce from the holomorphic Lefschetz formula that
$\tilde g$ has fixed points, hence also $g$.\cqfd
\medskip
\ind Recall that the moduli space $M_G$ is constructed as a {\it good quotient}
of a smooth scheme $R$ by a reductive group $\Gamma $ [Se] -- this implies in
particular that the closed points of $M_G$ correspond to the closed orbits of
$\Gamma $ in $R$. In order to compare the Picard groups of $\Pic(M_G)$ and
$\Pic({\cal M}_G)$, we will need a more precise result:

\smallskip
\th Lemma
\enonce There exists a presentation of
${\cal M}_{G}^{\delta,ss}$ as a quotient of a smooth scheme $R$ by a
reductive group $\Gamma $, such that the moduli space $M_{G}^\delta$ is
a good quotient of $R$ by $\Gamma$.
\endth\label{comppres}
{\it Proof}:   We will explain the proof in some detail
for the untwisted case,  then indicate how to adapt the argument  to
the general situation.
\ind We fix a faithful representation
$\rho:G\rightarrow {\bf SL}_r$ and an integer
$N$ such that for every semi-stable $G$\tx bundle $P$, the vector bundle
$P_{\rho}(Np)$ is generated by its global
sections and
satisfies $\ H^{1}(X,P _{\rho}(Np))=0$.
Let $M=r(N+1-g)$.  For any complex scheme
 $S$, we denote by $\underline{R}_G(S)$ the set of isomorphism classes of pairs
$(P ,\alpha)$, where $P $ is a   $G$\tx bundle on
$X\times S$ whose restriction to $X\times\{s\}$, for each closed point $s\in
S$, is
semi-stable, and
$\alpha:{\cal{O}}_{S}^M\iso pr_{1*}P _{\rho}(Np)$   an isomorphism. We
define in this way a functor $\underline{R}_G$ from the category of
${\bf C}$\tx schemes to the category of sets; we claim that it is
representable by a scheme $R_G$. If  $G={\bf SL}_r$, this follows from
Grothendieck theory of the Hilbert scheme [G1].  In the general case, we
observe that  reductions to $G$ of the structure group of a ${\bf SL}_r$\tx
bundle $P$ correspond canonically to global sections of the bundle $P/G$;
it follows that $\underline{R}_G$ is isomorphic to the functor of global
sections  of ${\cal P}/G$, where ${\cal P}$ is the  universal ${\bf
SL}_r$\tx bundle on  $X\times R_{{\bf SL}_r}$.  Again by [G1], this functor
is representable  by a scheme $R_G$, which is affine over $R_{{\bf SL}_r}$.

\ind Put $\Gamma ={\bf GL}_M$. The group  $\Gamma $ acts on ${\underline
R}_G$, and therefore on  $R_G$, by the rule $g\cdot(P ,\alpha)=(P ,\alpha
g^{-1} )$. This action lifts to
 the universal $G$\tx bundle
${\cal{P}}$ over $X\times R_G$ as follows: by construction  the
universal pair
$({\cal{P}},\alpha)$ is isomorphic to
$((\Id\times g)^{*}{\cal{P}},\alpha\rond g)$, hence there is an
isomorphism of $G$\tx bundles
$\sigma_{g}:(\Id\times g)^{*}{\cal{P}}\rightarrow {\cal P}$ such that
$\alpha \rond g^{-1}=pr_{1*}(\sigma_{g,\rho})\rond \alpha$. Since $\rho$ is
faithful this isomorphism is uniquely
determined by $pr_{1*}(\sigma_{g,\rho})$, hence depends
only on
$g$ and defines the required  lifting.

\ind The $\Gamma $\tx equivariant morphism
$\varphi_{{\cal{P }}}:R_G\rightarrow {\cal{M}}_G$
induces a morphism of stacks
$\overline{\varphi}_{{\cal{P}}}:[R_G/\Gamma ]\rightarrow {\cal{M}}_{G}^{ss}$
which is easily seen to be an isomorphism. We also have a $\Gamma $-equivariant
morphism $\psi_{\cal P} :R_G\rightarrow M_G$; if there exists
a good quotient $R_G//\Gamma $, the  universal property of $M_G$ implies
that $\psi _{\cal P}$ must induce an isomorphism of this quotient onto
$M_G$. The existence of such a good quotient is classical in the case
$G={\bf SL}_r$ (possibly after increasing $N$); for general $G$, since the
canonical map $R_G\rightarrow  R_{{\bf SL}_r}$ is $\Gamma $\tx equivariant
and affine, the existence of a good quotient of $R_{{\bf SL}_r}$ by $\Gamma
$ implies the same property for $R_G$ thanks to  a lemma of Ramanathan ([R1],
lemma
4.1).

\ind Let us finally consider the twisted case. We choose an embedding of
 the center $Z$ of $G$ in a torus $T={\bf G}_m^s$, and an embedding $\rho
:G\rightarrow \pprod_{i=1}^s {\bf GL}_{r_i}$ such that $\rho (Z)$ is central;
we put
$S=(\pprod_{}^{}  {\bf GL}_{r_i})\times (T/Z)$. The map  $(g,t)\mapsto
(t^{-1}\rho
(g), t\ \hbox{mod.}\ Z)$ of $G\times T$ into $S$
 defines an embedding of $C_ZG$ into $S$, which maps the center of $C_ZG
$ into the center of $S$, so that a $C_ZG$\tx bundle $P$ on $X$ is semi-stable
if and only if the associated $S$\tx bundle is semi-stable. We then argue as
before,
replacing ${\bf SL}_r$ by $S$.\cqfd
\smallskip
\th Proposition
\enonce  Assume that the group $G$ is almost simple. The group $\Pic
(M_G^\delta)$ is infinite cyclic, and the homomorphism  $\pi^*:\Pic
(M_G^\delta)\rightarrow \Pic (M_{\widetilde{G}}^\delta)$ is injective.
\endth\label{cyc}

\ind The second assertion follows from Lemma \ref{unirat} b) and
(\ref{char}); it is therefore enough to prove the first one when $G$ is
simply connected.  The proof then is the same as in the untwisted case
([L-S] or [K-N]): since the stack ${\cal M}_G^\delta$ is smooth, the
restriction map $\Pic({\cal M}_G^\delta)\rightarrow  \Pic({\cal
M}_G^{\delta,ss})$ is surjective, hence by  Prop.\ \ref{Pic-Bun} the group
$\Pic({\cal M}_G^{\delta,ss})$ is cyclic; it remains to prove that the pull
back homomorphism $\Pic(M_G^\delta)\rightarrow \Pic({\cal M}_G^{\delta,ss})$
is  injective.  \ind We choose a presentation of ${\cal M}_G^{\delta,ss}$ as a
quotient of a smooth scheme $R$ by a  reductive group $\Gamma $, such that the
moduli space ${ M}_G^\delta$ is a good  quotient of $R$ by $\Gamma $ (lemma
\ref{comppres}); then line bundles on ${\cal M}_G^{\delta,ss}$ (resp.\ on
$M_G^\delta$) correspond to line bundles on $R$ with a $\Gamma $\tx
linearization (resp.\ a $\Gamma $\tx linearization $\sigma$ such that $\sigma
(\gamma )_r=1$ for each $(\gamma ,r)\in \Gamma \times R$ such that $\gamma
r=r$), hence our assertion.\cqfd
\medskip
\ind  In what follows we will identify the group $\Pic (M_G^\delta)$ with its
image in $\Pic (M_{\widetilde{G}}^\delta)$; our aim will be to find its
generator.

\vskip1cm
\section{The Picard groups of $M_{{\bf Spin}_{r}}$ and  $M_{G_2}$}
\label{Pic(M_Spin)}
\subsection  In this section we complete the results of [L-S] in the simply
connected case. The cases  $G={\bf SL}_r$ or ${\bf Sp}_{2l}$ are dealt with in
{\it loc. cit.}. We now consider the case $G={\bf Spin}_r$; we
denote by ${\cal D}$ the determinant bundle on $M_{{\bf Spin}_{r}}$ associated
to the standard representation $\sigma $ of ${\bf Spin}_{r}$ in ${\bf C}^r$.

\th Proposition
\enonce Let $r$ be an integer $\geq 7$. The group
$\Pic(M_{{\bf Spin}_r})$ is generated by ${\cal D}$.
\endth\label{Pic_Spin}

{\it Proof}: Choose a presentation of ${\cal M}_{{\bf Spin}_{r}}^{ss}$ as a
quotient of a smooth scheme $R$ by a reductive group $\Gamma $, such that
$M_{{\bf Spin}_r}$ is a good  quotient of $R$ by $\Gamma $   (lemma
\ref{comppres}). Let ${\cal S}$
be the universal  ${\bf Spin}_{r}$\tx  bundle  on
$X\times R$. We fix a theta-characteristic $\kappa $ on $X$; this
allows us to define the pfaffian line bundle ${\cal P}_\kappa $ on $R$,
which is a square root of $\det Rpr_{2*}({\cal S}_\sigma \otimes\kappa )$
[L-S]. The action of $\Gamma $ on ${\cal S}$ defines a $\Gamma $\tx
linearization of ${\cal P}_\kappa $.
\ind By [L-S] we know that the group
of $\Gamma $\tx linearized line bundles on $R$ (isomorphic to $\Pic({\cal
M}_{{\bf Spin}_{r}}^{ss})$) is generated by ${\cal P}_\kappa $,  so all we
have to prove is that ${\cal P}_\kappa $ itself does not descend to
$R//\Gamma $, i.e.\  to exhibit a  closed point $s\in R$ whose stabilizer
in $\Gamma $ acts nontrivially on the fibre of ${\cal P}_\kappa $ at $s$.
If $s$ corresponds to a semi-stable ${\bf Spin}_r$\tx bundle $P$, its
stabilizer is the group $\Aut(P)$;  since the formation of pfaffians
commutes with base change, its action on $({\cal P}_\kappa)_s $ is the
natural action of $\Aut(P)$ on $\bigl(\wedge^{\rm max}H^0(X,P_\sigma
\otimes\kappa )\bigr)^{-1}$ [L-S].  \ind To construct $P$ we follow [L-S],
prop.\ 9.5: we choose  a stable ${\bf SO}_4$\tx bundle $Q$ and a stable
${\bf SO}_{r-4}$\tx bundle $Q'$ with $w_2(Q)=w_2(Q')=1$.
Let $H$ be the subgroup ${\bf SO}_{4}\times {\bf SO}_{r-4}$ of ${\bf
SO}_{r}$, and $\widetilde{H}$ its inverse image in ${\bf Spin}_r$. The
$H$\tx bundle $Q\times Q'$ has $w_2=0$ by construction, hence admits a
$\widetilde{H}$\tx structure; we choose one, and take for $P$ the associated
${\bf Spin}_r$\tx bundle. Let $\gamma $ be  a central element  of
$\widetilde{H}$ lifting the element $(-1,1)$ of $H$. Then $\gamma $ defines
an automorphism of $P$, which acts on the associated vector bundle
$P_\sigma =Q_\sigma \oplus Q'_\sigma $ as $(-\Id,\Id)$ (we use the same
letter
$\sigma$ to denote the standard representation of all orthogonal group in
sight). Therefore $\gamma $ acts on $({\cal P}_\kappa)_s$ by multiplication
by
$(-1)^{h^0(Q_\sigma \otimes\kappa )}$. But $h^0(Q_\sigma \otimes\kappa )$ is
congruent to $w_2(Q)$ (mod.$\,2$) [L-S, 7.10.1], hence our assertion.\cqfd

\medskip
\rem{Remark}  For $r\le6$ the group ${\bf Spin}_r$ is of type A or C, so we
already have a complete description of $\Pic(M_{{\bf Spin}_r})$. It is worth
noticing that the above result does not hold for $r\le 6$: using the
exceptional isomorphisms  one checks easily that  $\Pic(M_{{\bf Spin}_r})$
is generated by a square root of ${\cal D}$ for $r=5$ or $6$ and a fourth
root for $r=3$ -- while it is isomorphic to ${\bf Z}^2$ for $r=4$.

\medskip
\ind We now consider the case when $G$ is of type $G_2$. The group $G$ is
the automorphism group of the octonion algebra ${\bf O}$ over ${\bf C}$
([Bo], Alg\`ebre III, App.); in particular it
has a natural orthogonal representation $\sigma $ in the $7$\tx dimensional
space
${\bf O}/{\bf C}$. We denote by ${\cal D}$ the
determinant bundle on $M_G$ associated to this representation.

\th Proposition
\enonce The group $\Pic(M_G)$ is generated by ${\cal D}$.
\endth\label{G_2}
{\it Proof}: As before we choose a presentation of ${\cal M}_{G}^{ss}$
as a quotient of a smooth scheme $R$ by a reductive group $\Gamma $, such that
$M_{G}=R//\Gamma $;
choosing a theta-characteristic $\kappa $ on $X$
allows to define a pfaffian line bundle ${\cal P}_\kappa $ on
$R$, with a natural  $\Gamma $\tx linearization. By [L-S], thm.\ 1.1,
${\cal P}_\kappa $ generates the group of $\Gamma $\tx linearized line
bundles on $R$;  we must again prove that it does not
descend to $R//\Gamma $, i.e.\ exhibit a $G$\tx bundle $P$ such that
$\Aut(P)$ acts nontrivially on $\wedge^{\rm max}H^0(P_\sigma \otimes\kappa
)$.

\ind Let  $V$ be a
$3$\tx dimensional vector space over ${\bf F}_2$. The algebra ${\bf O}$ has a
basis
$(e_\alpha )_{\alpha \in V}$, with multiplication rule $$e_\alpha \,e_\beta
=\varepsilon (\alpha ,\beta )\,e_{\alpha +\beta }\ ,$$for a certain function
$\varepsilon :V\times V\rightarrow \{\pm 1\}$. Suppose given a homomorphism
$\alpha
\mapsto L_\alpha $ of $V$ into $J$. We view $J$ as the moduli space for degree
$0$
line bundles with a trivialization at $p$; for each pair $(\alpha ,\beta )$ in
$V$
we have a unique isomorphism $u_{\alpha \beta }:L_\alpha \otimes L_\beta
\rightarrow
L_{\alpha +\beta }$ compatible with these trivializations. We endow the
${\cal O}_X$\tx module ${\cal A}=\sdir_{\alpha \in V}^{}L_\alpha $ with the
algebra
structure defined by the map ${\cal A}\otimes{\cal A}\rightarrow {\cal A}$
which
coincides with $\varepsilon (\alpha ,\beta )\,u_{\alpha \beta }$ on $L_\alpha
\otimes L_\beta $. It is a sheaf of ${\cal O}_X$\tx algebras, locally
isomorphic to
${\cal O}_X\otimes_{\bf C}{\bf O}$. Let $P$ be the associated $G$\tx bundle
(the
sections of $P$ over an open subset $U$ of $X$ are algebra isomorphisms of
${\cal
O}_U\otimes_{\bf C}{\bf O}$ onto ${\cal A}_{|U}$). Since the pull back of $P$
to any
finite covering of $X$ on which the $L_\alpha $'s are trivial is trivial, $P$
is
semi-stable. The vector bundle $P_\sigma $ is simply $\sdir_{\alpha
\not=0}^{}L_\alpha $. Let $\chi :V\rightarrow \{\pm 1\}$ be a nontrivial
character;
the diagonal endomorphism $(\chi (\alpha ))_{\alpha \in V}$ of ${\cal A}$ is an
algebra automorphism, and therefore defines an automorphism $\iota $ of $P$,
which
acts on $P_\sigma $  with eigenvalues $(\chi (\alpha ))_{\alpha
\not=0}$. Hence $\iota $ acts on $\wedge^{\rm max}H^0(P_\sigma \otimes\kappa )$
by
multiplication by $(-1)^h$, with $h=\sum_{\chi (\alpha )=-1} h^0(L_\alpha
\otimes\kappa )$. Since the function $\alpha \mapsto h^0(L_\alpha
\otimes\kappa )$ (mod.$\,2$) is quadratic, an easy computation gives that $h$
is even
if and only if the image of $\Ker \chi $ in $J_2$ is totally isotropic with
respect
to the Weil pairing.
 Clearly we can choose our map $V\rightarrow  J_2$ and the
character $\chi $ so that this does not hold; this provides the required
example.\cqfd
\vskip1cm
\section
{The Picard group of $M^0_G$}

\ind In the study of $\Pic(M_G^\delta)$, contrary to
what we found for the moduli stacks, the degree $\delta$ plays an important
role. We treat first the degree $0$ case, which is easier. Let us start with
the case $A_l$. We recall that the determinant bundle ${\cal D}$ exists on
the moduli space $M_{{\bf SL}_r}$, and generates its Picard group.

\smallskip
\th Proposition
\enonce Let $G={\bf SL}_r /\Mu_s$, with $s$ dividing $r$.

{\rm a)} If $s$ is odd or $r/s$ is even,  $\Pic (M_G^0)$ is generated by ${\cal
D}^s$.

{\rm b)} If $s$ is even and $r/s$ is odd, $\Pic (M_G^0)$ is generated by ${\cal
D}^{2s}$.

\ind  In particular, $\Pic( M_{{\bf PGL}_r}^0)$ is generated by ${\cal D}^{r}$
if $r$ is odd and by ${\cal D}^{2r}$ if $r$ is even.

\endth\label{A0}
{\it Proof}: We identify $M_{{\bf SL}_r }$ with the moduli space of
semi-stable vector bundles of rank $r$ and trivial determinant on $X$. Let
$J^{g-1}$ be the component of the Picard variety of $X$ parameterizing line
bundles of degree $g-1$, and $\Theta \i J^{g-1}$ the canonical theta divisor.
It is shown in [B-N-R] that for $E$ general in $M_{{\bf SL}_r }$, the
 condition $H^0(X,E\otimes L)\not=0$ defines a divisor $D(E)$ in $J^{g-1}$
which belongs to the linear system $|r\Theta |$, and that the rational map
$D:M_{{\bf SL}_r }\dasharrow |r\Theta |$ thus defined satisfies $D^*{\cal
O}(1)={\cal D}$. Using (\ref{list} {\it a})  we deduce that the alternate
form associated to ${\cal E}(J_r,{\cal D})$ is the inverse of the Weil pairing
$e_r$;  its restriction to $J_s$ is $e_s^{-r/s}$ (\ref{restric}).
 From Prop. \ref{desc-prop}, we conclude that the line bundles ${\cal
D}^s$ in case a) and ${\cal D}^{2s}$ in case b) descend to $M_G^0$.
\ind It remains to prove that these line bundles are indeed in each case
generators of $\Pic (M_G^0)$. Consider first the case $s=r$. Since  the
extension  ${\cal E}(J_r,{\cal D})$ is of order $r$ in $H^2(J_r,{\bf
C}^*)$, the smallest power of ${\cal D}$ which admits a $J_r$\tx
linearization is ${\cal D}^r$, so the conclusion follows from Prop.\
\ref{desc-prop}. In the general case, put $M:=M_{{\bf SL}_r }$, and assume
that some power ${\cal D}^k$ of ${\cal D}$ descends to  $M/J_s$. Observe that
$M/J_r$ can be viewed as the quotient of $M/J_s$ by $J_{r/s}$.

\ind Assume that $r/s$ is even. We know by Prop.\ \ref{desc-prop} that ${\cal
D}^{2kr/s}$ descends to $M/J_r$; since $r$ is even, this implies by what we
have
seen that $2r$ divides $2kr/s$, hence that $k$ is a multiple of $s$. If $r/s$
is odd, then ${\cal D}^{kr/s}$ descends by Prop.\ \ref{desc-prop}, and
therefore
$k$ is a multiple of $s$ or $2s$ according to the parity of $r$.\cqfd

\medskip
\subsection We now consider the case of the orthogonal and symplectic group.
If $G={\bf SO}_r$ or ${\bf Sp}_r$ ($r$ even), we will denote by ${\cal D}$ the
determinant bundle on $M_G$, i.e.\ the pull back of the determinant bundle on
$M_{{\bf SL}_r}$ by the morphism associated to the standard representation. We
know
that the group $\Pic(M_{{\bf Sp}_{r}})$ is generated by ${\cal D}$ ([L-S],
1.6), and
 that $\Pic(M_{{\bf Spin}_{r}})$ is generated by the pull back of ${\cal D}$
(Prop.\ \ref{Pic_Spin}); it follows that the Picard group of each component of
$M_{{\bf SO}_r}$ is generated  by ${\cal D}$. It remains to consider the groups
${\bf
PSp}_{2l} $ and ${\bf PSO}_{2l} $.

\th Proposition
\enonce Let $G={\bf PSp}_{2l} $
or ${\bf PSO}_{2l}\quad  (l\geq 2)$.

{\rm a)}  If $l$ is  even,  $\Pic (M_G^0)$ is generated by ${\cal D}^2$.

{\rm b)}  If $l$ is odd,  $\Pic (M_G^0)$ is generated by ${\cal D}^4$.
\endth
\label{BC0}
{\it Proof}: The extension ${\cal E}(J_2,{\cal D})$ is the pull back to $J_2$
of
the Heisenberg extension of $J_{2l}$, and  the corresponding alternate form is
$e_2^l$ (\ref{restric}). We deduce from Proposition \ref{desc-prop} that ${\cal
D}^{2}$ descends to  $M_G^0$ if $l$ is even, and that ${\cal D}^4$ descends but
${\cal D}^2$ does not if $l$ is odd.
\ind It remains to prove that ${\cal D}$ does not descend when $l$ is even.
Let us consider for instance the case of the symplectic group; for every
integer
$n$, we put $M_n=M_{{\rm\bf Sp}_{2n}}$ and denote by ${\cal D}_n$ the
determinant
line bundle on $M_n$. Write $l=p+q$, where $p$ and $q$ are odd (e.g.\ $p=1$,
$q=l-1)$, and consider the morphism $u:M_p\times M_q\rightarrow
M_l$ given by $u((E,\varphi ),(E',\varphi '))=(E\oplus E',\varphi
\oplus \varphi ')$. It is $J_2$\tx equivariant and satisfies $u^*{\cal
D}_l={\cal
D}_p\boxtimes{\cal D}_q$.
The group $J_2\times J_2$ acts on $M_p\times M_q$; from (\ref{list} {\it c})
one
deduces that the alternate form $e$ corresponding to the  extension ${\cal
E}(J_2\times J_2,{\cal D}_p\boxtimes{\cal D}_q)$ is given by $e((\alpha ,\alpha
'),(\beta ,\beta '))=e_2(\alpha ,\beta )\,e_2(\alpha ',\beta ')$.
If ${\cal D}_l$ descends to $M_l/J_2$, then ${\cal D}_p\boxtimes{\cal D}_q$
descends to $(M_p\times M_q)/J_2$, and we can apply (\ref{list} {\it d}) to
the variety $M_p\times M_q$ and the diagonal embedding $J_2\i J_2\times
J_2$. We conclude that the form $e$ is the pull back of an alternate form on
$J_2$ by the sum map $J_2\times J_2\rightarrow J_2$.  This is clearly
impossible, which proves that  ${\cal D}_l$ does not descend to $M_l/J_2$.
\ind The same argument applies to the orthogonal groups, except that one
has to be  careful about the definition of $M_1$: we take it to be the
Jacobian of $X$, by associating to a line bundle $\alpha $ on $X$ the vector
bundle $\alpha \oplus\alpha ^{-1}$ with the standard isotropic form. Then
${\cal D}_1$ is the line bundle ${\cal O}(2\Theta )$. The alternate form
associated to ${\cal E}(J_2,{\cal D}_1)$ is $e_2$, and the rest of the
argument applies without any change.\cqfd

\medskip
\rem{Remark}
 There remains one case to deal with.
 When $l$ is
even, the center $Z$ of ${\rm\bf Spin}_{2l}$ is isomorphic to
$\Mu_2\times\Mu_2$, so
it contains two subgroups of order $2$ (besides the kernel of the homomorphism
${\bf Spin}_{2l}\rightarrow {\bf SO}_{2l}$). These subgroups  are
exchanged by the outer automorphisms of ${\bf Spin}_{2l}$, so the
corresponding  quotient groups are canonically
isomorphic; let us denote them by $G$. Since $M_G^0$  dominates $M^0_{{\bf
PSO}_{2l}}$, it follows from Prop.\ \ref{BC0} that ${\cal D}^2$ descends to
$M_G^0$. If $l$ is not divisible by $4$, one can show that ${\cal D}$ does
not descend to $M_G^0$, so $\Pic(M_G^0)$ {\it is generated by} ${\cal
D}^2$. If $l=4$, one sees using the triality automorphism that ${\cal D}$
descends; we do not know what happens for $l=4m$, $m\ge 2$.

 \vskip1cm
\section{The Picard group of $M_{{\bf PGL}_r}^d$}

\ind In this section we consider the  component $M_{{\bf PGL}_r}^d$  of the
moduli space $M_{{\bf PGL}_r}$, for $0<d<r$.  It is
the quotient by  $J_r$ of the moduli space $M_{{\bf SL}_r}^d$ of semi-stable
vector
bundles of rank $r$ and  determinant ${\cal O}_X(dp)$. We denote by $\delta$
the
g.c.d. of $r$ and $d$. If $A$ is a  vector bundle on $X$ of rank $r/\delta$
and degree $(r(g-1)-d)/\delta$ which is general enough, the condition
$H^0(X,E\otimes
A)\not=0$ defines a Cartier divisor on $M_{{\bf SL}_r}^d$; the associated line
bundle
${\cal L}$  (sometimes called the {\it theta line bundle}) is independent of
the
choice of $A$, and generates $\Pic(M_{{\bf SL}_r}^d)$ [D-N].

\medskip
\th Proposition
\enonce  The group $\Pic(M_{{\bf PGL}_r}^d)$ is generated
by ${\cal L}^\delta$ if $r$ is odd and  by ${\cal L}^{2\delta}$ if $r$ is even.
\endth
\label{A_d}
\ind Choose a stable vector bundle $A$ of rank $r/\delta$
and determinant ${\cal O}_X(-{d\over \delta}p)$, and consider the morphism
$a:E\mapsto E\otimes A$ of $M_{{\bf SL}_r}^d$ into $M_{{\bf
SL}_{r^2/\delta}}^0$. By
definition ${\cal L}$ is the pull back of the determinant bundle ${\cal D}$ on
the
target. The map $a$ is $J_r$ equivariant, hence induces an isomorphism
${\cal E}(J_r,{\cal D})\cong {\cal E}(J_r,{\cal L})$. We have seen in the proof
of
Prop.\ \ref{A0} that the  alternate form associated to ${\cal E}(J_r,{\cal D})$
is
$e_r^{-r/\delta}$; hence the smallest power of
 ${\cal  L}$ which  descends to
$M_{{\bf PGL}_r}^d$ is ${\cal L}^\delta$.
Therefore  it is enough to prove that ${\cal
L}^\delta$ descends to  $M_{{\bf PGL}_r}^d$ when $r$ is odd and that ${\cal
L}^{2\delta}$ but not ${\cal  L}^{\delta}$ descends when $r$ is even.

\ind We will prove this by reducing to the degree $0$ case with the
help of the Hecke correspondence.  Let us denote simply by $M$ the moduli space
$M_{{\bf SL}_r}^1$ of stable vector bundles of rank $r$ and determinant ${\cal
O}_X(p)$ on $X$. There exists a Poincar\'e bundle  ${\cal  E}$ on $X\times M$;
we
denote by ${\cal E}_p$ its restriction  to $\{p\}\times M$, viewed as a vector
bundle on $M$.  We fix  an integer $h$ with  $0<h<r$ and
let
 ${\cal P}={\bf G}_M(h,{\cal
E}_p)$  the Grassmann bundle parameterizing rank $h$ locally free quotients of
${\cal E}_p$. A point of ${\cal P}$ can be viewed as a pair
$(E,F)$ of vector bundles with $E\in M$, $E(-p)\i F\i E$ and $\dim
(E_p/F_p)=h$.
\smallskip

\th Lemma
\enonce If  $E$ is general enough in $M$, for any pair
$(E,F)$ in ${\cal P}$ the vector bundle $F$ is semi-stable, and stable if
$g\ge3$.
\label{slopes}
\endth
{\it Proof}: We will actually prove  a more precise result. If $G$ is a
vector bundle on $X$, define the {\it stability degree} $s(G)$ of $G$ as
the minimum of the rational numbers $\mu(G'')-\mu(G')$ over all exact
sequences $0\rightarrow G'\rightarrow G\rightarrow G''\rightarrow 0$. One
has $s(G)\ge0$ if and only if $G$ is semi-stable, $s(G)>0$ if and only if
$G$ is stable, and $s(G)=g-1$ when $G$ is a general stable vector bundle
[L, Hi].  \ind Let $E,F$ be two vector bundles on $X$, with $E(-p)\i F\i
E$. The lemma will follow from the inequality $$s(F)\ge s(E)-1$$(note that
since $E$ and $F$ play a
symmetric role, this implies $|s(E)-s(F)|\le1$). Let $Q_p$ be the sheaf
$E/F$ (with support $\{p\}$), and $h$ the dimension of its fibre at $p$.
 Let $F'$ be a subbundle of $F$, of rank $r'$. From the exact sequence
$0\rightarrow F/F'\rightarrow E/F'\rightarrow Q_p\rightarrow 0$ we get
  $$\mu(F/F')-\mu(F')=\mu(E/F')-{h\over r-r'}-\mu(F')\ge s(E)-{h\over
r-r'}\ .$$ Let $K_p:=\Ker(E_p\rightarrow Q_p)$. The exact sequence
$0\rightarrow E(-p)\rightarrow F\rightarrow K_p\rightarrow 0$ induces  an
exact sequence $0\rightarrow E'\rightarrow F'\rightarrow K_p$, with
$E':=F'\cap E(-p)$. Therefore $$\mu(F/F')-\mu(F')\ge\mu
(E(-p)/E')-\mu(E')-{r-h\over r'}\ge s(E)-{r-h\over r'}\ .$$ Since  one of
the two numbers $\displaystyle {h\over r-r'}$ and $\displaystyle {r-h\over
r'}$ is $\le1$, we get the required inequality.\cqfd

\medskip
\ind Let us denote by $M'$ the moduli space $M_{{\bf SL}_r}^{1-h}$. Using
the lemma we get a diagram
\input xypic $$\diagram
 & {\cal  P}\dlto_{q} \xdashed[dr]^{q'}|>\tip &\\
  M &&M'\\
\enddiagram$$
\smallskip
(``Hecke diagram"), where $q$ (resp.\ $q'$) associates to a  pair $(E,F)$
the vector bundle $E$ (resp.\ $F$, provided $F$ is semi-stable).
 \ind Let ${\cal L}$ and ${\cal L}'$ be the theta line bundles on
$M$ and $M'$.  Let $\delta$ be the g.c.d. of
$r$ and $1-h$.

\smallskip
\th Lemma
\enonce One has $K_{\cal  P}=q^*{\cal  L}^{-1}\otimes
q'^*{\cal  L}'^{-\delta}$.
\label{KP}
\endth
{\it Proof}: Let $E$ be a general vector bundle in $M$; let us  compute
the restriction of $q'^*{\cal L}'$ to the fibre $q^{-1}(E)$.   On
$X\times{\cal P}$ we have a canonical exact sequence $$0\rightarrow {\cal
F}\longrightarrow (1_X\times q)^*{\cal E}\longrightarrow (i_p)_*{\cal
Q}_p\rightarrow 0\ ,$$ where ${\cal Q}_p$ is the universal quotient bundle of
$q^*{\cal E}_p$ on ${\cal P}$ and $i_p$ the embedding of ${\cal P}=
\{p\}\times{\cal P}$ in $X\times{\cal P}$. For each point $P=(E,F)$ of ${\cal
P}$
 this exact sequence gives by  restriction to $X\times\{P\}$  the exact
sequence
$0\rightarrow F\rightarrow E\rightarrow Q_p\rightarrow 0$ defining $P$; in
particular, one has   ${\cal F}_{X\times\{P\}}=F$, and the map $q':{\cal
P}\dasharrow M'$ is the classifying map associated to ${\cal F}$. It
follows that $q'^*{\cal L}'$ is the determinant bundle associated to ${\cal
F}\otimes A$,  where $A$ is a vector bundle
of rank $r/\delta$ and appropriate degree.
\ind Now let $E\in M$; put ${\bf G}=q^{-1}(E)={\bf G}(h,E_p)$, and denote by
 $\pi:X\times{\bf G}\rightarrow X$ and $\rho :X\times{\bf G}\rightarrow
{\bf G}$  the two projections. The
restriction of the above exact sequence to $X\times{\bf G}$ gives, after tensor
product with $\pi^*A$, an exact sequence$$0\rightarrow {\cal
F}\otimes\pi^*A\longrightarrow \pi^*({\cal E}\otimes A)\longrightarrow
(i_p)_*{\cal
Q}_p^{r/\delta}\rightarrow 0\ ;$$
applying
$R\rho _*$  and taking
determinants, we obtain
$$\det R\rho _*({\cal F}\otimes \pi^*A)\cong (\det {\cal
Q}_p)^{r/\delta}={\cal O}_{\bf G}(r/\delta)\ .$$
\ind The restriction of $K_{\cal P}$ to ${\bf G}$ is $K_{\bf G}={\cal
O}_{\bf G}(-r)$; since $\Pic({\cal P})$ is generated by ${\cal O}_{\cal P}(1)$
and $q^*\Pic(M)$, one can write $K_{\cal P}=q'^*{\cal L}'^{-\delta}\otimes
q^*{\cal L}^a$ for some integer $a$. To compute $a$ we consider  the
restriction of $q^*{\cal L}$ to a general fibre $q'^{-1}(F)$: by lemma
(\ref{slopes}) this fibre can be identified with the Grassmann variety ${\bf
G}(r-h,F_p)$, and the same argument as above shows that the restriction of
$q^*{\cal L}$ is equal to ${\cal O}_{\bf G}(r)$, that is to the restriction of
$K_{\cal P}^{-1}$. This gives $a=-1$, hence the lemma.\cqfd

 \medskip \ind Observe that the group $J_r$
acts in a natural way on ${\cal P}$, by the rule $\alpha \cdot
(E,F)=$ $(E\otimes\alpha ,F\otimes\alpha )$; the Hecke diagram  is $J_r$\tx
equivariant.

\smallskip
\th Lemma
\enonce Let $s$ be an integer dividing $r$. The canonical bundle $K_{\cal P}$
descends to ${\cal P}/J_s$, except if $s$ is even and $h$ and $r/s$ are odd; in
this last case $K_{\cal P}^2$ descends.
\endth
\label{desc-KP}
{\it Proof}: {\it a}) We first prove that $K_M$ descends to $M/J_r$. Let $\pi$
and $\rho $ denote the projections from $X\times M$ onto $X$ and $M$
respectively.
By  deformation theory, the tangent bundle $T_M$ is canonically isomorphic
to $R^1\rho _*({\cal E}nd_0({\cal E}))$, where ${\cal E}nd_0$ denote the
sheaf of traceless endomorphisms; it follows that $K_M$ is the inverse of
the determinant bundle $\det R\rho _*({\cal E}nd_0({\cal E}))$. Since
${\cal E}nd_0({\cal E})$ has trivial determinant, this is also equal to
$\det R\rho _*({\cal E}nd_0({\cal E})\otimes \pi^*L)$ for any line bundle
$L$ on $X$ (see e.g.\ [B-L], 3.8); therefore $K_M^{-1}$ is the pull back of
the generator ${\cal L}$ of $\Pic(M_{{\bf SL}_{r^2-1}})$   by the morphism
$M\rightarrow M_{{\bf SL}_{r^2-1}}$ which maps $E$ to ${\cal
E}nd_0(E)$. This morphism  factors through the quotient $M/J_r$, hence our
assertion.

 \ind {\it b}) Therefore we need only to consider the relative canonical bundle
$K_{{\cal P}/M}$, with its  canonical $J_r$\tx linearization. Let $\alpha \in
J_r$, and let $P=(E,F)$ be a fixed point of $\alpha $ in ${\cal P}$; we want
to compute the tangent map $T_P(\alpha)$ to $\alpha$ at $P$. The vector bundle
$E\in M$ is  fixed by $\alpha $, and the action of $\alpha $ on the
fibre
$q^{-1}(E)={\bf G}(E_p)$ is induced by the automorphism $\tilde \alpha $ of
$E_p$
obtained from the  isomorphism $\varphi_\alpha  :E\rightarrow E\otimes\alpha $
(note that $\varphi_\alpha  $, hence also $\tilde \alpha $, are uniquely
determined up to a scalar, since $E$ is stable).
\ind Let $0\rightarrow K_p\rightarrow
E_p\rightarrow Q_p\rightarrow 0$ be the exact sequence corresponding to $P$.
The
tangent space to ${\bf G}(E_p)$ at $P$ is canonically isomorphic to
$\Hom(K_p,Q_p)$, hence its determinant is canonically isomorphic to
$(\det E_p)^{-h}(\det Q_p)^{r}$;  we conclude that $\det T_P(\alpha )$  is
equal
to
$(\det\tilde
\alpha)^h
$, where $\tilde\alpha$ is normalized so that $\tilde\alpha^r=1$.

\ind {\it c}) It remains to compute $\det\tilde \alpha$. Now the fixed points
of
$\alpha $ on $M$ are easy to describe [N-R]: let $s$ be the order of $\alpha $,
and $\pi:\widetilde{X}\rightarrow X$ the associated  \'etale $s$\tx sheeted
covering;
 a vector bundle $E$ on $X$ satisfies $E\otimes\alpha \cong E$ if and only
if it is of the form $\pi_*\widetilde{E}$ for some vector bundle
$\widetilde{E}$ on
$\widetilde{X}$, of rank $r/s$.  To evaluate $\varphi _\alpha $ at  $p$, we can
trivialize $\widetilde{E}$ in a neighborhood of $\pi^{-1}(p)$: write
$\widetilde{E}=\pi^*T$, where $T={\cal O}_X^{r/s}$. Then one has
$\pi_*\widetilde{E}=\sdir_{i\in{\bf Z}/s}^{} T\otimes\alpha ^i$, and the
isomorphism $\varphi _\alpha $ maps identically $T\otimes\alpha ^i$ onto
$(T\otimes\alpha ^{i-1})\otimes\alpha $. It follows that  the
eigenvalues of $\tilde \alpha $  are the $s$\tx th roots of $1$, each
counted with multiplicity $r/s$. This implies in particular
 $\det \tilde \alpha =\zeta ^{r(s-1)/2}$, where  $\zeta $
is a primitive $s$\tx th root of $1$, and therefore $\det T_P(\alpha
)=(-1)^{h(s-1){r\over s}}$. The lemma follows.\cqfd

\medskip
\subsec{\it Proof of Proposition} \ref{A_d}:

\ind We first observe that a line bundle $L$ on $M$  descends to
$M/J_s$ if and only if its pull back to ${\cal P}$ descends to
${\cal P}/J_s$. In fact, we know by (\ref{list} {\it a}) that  $G$\tx
linearizations of $L$ correspond bijectively by pull back to  $G$\tx
linearizations of $q^*L$; for $\alpha \in J_s$, any fixed point $E$ of
$\alpha $ in $M$ is the image of a point $P\in{\cal P}$ fixed by
$\alpha $, so with the notation of (\ref{char}) one has $\chi _E(\alpha )=\chi
_P(\alpha )$, which implies our assertion.
\ind Similarly,  a line bundle on $M'$  descends to $M'/J_s$
if and only if its pull back to ${\cal P}$ descends to ${\cal P}/J_s$:
what we have to check in order to apply the same argument is that every
component of the fixed locus ${\rm Fix}_{M'}(\alpha)$ is dominated
by a component of ${\rm Fix}_{{\cal P}}(\alpha)$, and conversely that
every component of ${\rm Fix}_{{\cal P}}(\alpha)$ dominates a component
of ${\rm Fix}_{M'}(\alpha)$. But this follows easily from the
description of the fixed points of $\alpha $ given above (\ref{desc-KP} {\it
c}).

\ind We first consider the case $h=1$. If $r$ is odd, we know from Prop.\
\ref{A0}  and lemma \ref{desc-KP}  that ${\cal  L}'^r$ and $K_{{\cal
P}}=q^*{\cal  L}^{-1} \otimes q'^*{\cal L}'^{-r}$ descend to ${\cal P}/J_r$; it
follows that ${\cal L}$ descends to $M/J_r$. Assume that $r$ is even. Endow
$K_{\cal P}$ with its canonical $J_r$\tx linearization, ${\cal L}^r$ with the
$J_r$\tx linearization defined in (\ref{desc-prop}), and  $q^*{\cal L}$ with
the
$J_r$\tx linearization deduced from the isomorphism $K_{{\cal
P}}\cong q^*{\cal  L}^{-1} \otimes q'^*{\cal L}'^{-r}$. Let
$\alpha $ be  an
element of order $r$ in $J_r$, and $P$ a fixed point of $\alpha $ in ${\cal
P}$; we know that $\alpha $ acts on $(K_{\cal P})_P$ by
multiplication by $-1$ (\ref{desc-KP} {\it c}) and on $(q'^*{\cal L}'^r)_P$ by
multiplication by $\varepsilon (\alpha )$ (\ref{desc-prop}), hence it acts on
$(q^*{\cal L})_P$ by multiplication by $-\varepsilon (\alpha )$. Since
$-\varepsilon (\alpha +\beta )\not=(-\varepsilon (\alpha )\,(-\varepsilon
(\beta )))$
when $\alpha $ and $\beta $ are two elements of order $r$ orthogonal for the
Weil
pairing, we conclude that ${\cal L}$ does not descend, while of course ${\cal
L}^2$ descends. \ind We now apply the same argument with $h$ arbitrary. If $r$
is
odd, $K_{\cal P}$ and $q^*{\cal L}$ descend, hence ${\cal L}'^{\delta}$
descends. If $r$ is even, we get a $J_r$\tx linearization on $q'^*{\cal
L}'^\delta$
such that an element $\alpha$ of order $r$ in $J_r$ acts by multiplication by
$(-1)^{h+1}\varepsilon (\alpha )$; again this implies that ${\cal
L}'^{\delta}$ does not descend, while ${\cal L}'^{2\delta}$ descends.\cqfd
\bigskip
\rem{Remark}
The methods of this section allow to treat more generally in most cases  the
group ${\bf SL}_r/\Mu_s$, for $s$ dividing $r$. We will contend ourselves with
an
example, which we will need below:  the case $G={\bf SL}_{2l}/\Mu_2$ ($l\ge
1$).
The moduli space $M_G$ has two components, namely $M^0_G$ (treated in Prop.\
\ref{A0}) and the quotient $M_G^l$ of $M_{{\bf SL}_{2l}}^l$ by $J_2$.
The theta line bundle ${\cal L}$ on $M_{{\bf SL}_{2l}}^l$  is the pull back of
the
determinant bundle on $M_{{\bf SL}_{4l}}^0$ under the map $E\mapsto E\otimes
A$,
where $A$ is a stable vector bundle of  rank $2$ and degree $-1$. It follows
from
Prop.\ \ref{A0} that ${\cal L}^2$ descends  to $M_G^l$; on the other hand, by
Prop.\
\ref{A_d},  ${\cal L}^l$ and therefore ${\cal L}$ do not descend if $l$ is odd.
We
shall now prove that  ${\cal L}$ {\it descends to $M_G^l$ when $l$ is even}.

\ind Let $\lambda :M_{{\bf SL}_{2l}}^l\rightarrow M_{{\bf
SL}_{l(2l-1))}}$ be the morphism $E\mapsto\wedge^2E(-p)$. One checks easily
 that the pull back of the determinant bundle ${\cal D}$ on $M_{{\bf
SL}(l(2l-1))}$ is  ${\cal L}^{l-1}$ (e.g.\ by pulling back to the moduli stack,
and
using the fact that the Dynkin index of the representation $\wedge^2$ is
$2l-2$).
Now $\lambda $ factors through $M_{{\bf SL}_{2l}}^l/J_2$, therefore ${\cal
L}^{l-1}$
descends to this quotient. When $l$ is even, this implies that ${\cal L}$
itself
descends. \label{M(2l,l)}

\vskip1cm
\section {The Picard groups of $M_{{\bf PSp}_{2l}}$ and $M_{{\bf PSO}_{2l}}$}
\subsection In the case $C_l$, it remains only to consider the
component $M_{{\bf PSp}_{2l}}^1$, which is the quotient by $J_2$ of the
moduli space $M_{{\bf PSp}_{2l}}^1$  of semi-stable pairs $(E,\varphi )$,
where $E$ is a vector bundle of rank $2l$ on $X$ and  $\varphi
:\wedge^2E\rightarrow {\cal O}_X(p)$  a non-degenerate alternate form.
   Let ${\cal L}$ denote the pull back of the theta line bundle  by
the  natural map $M_{{\bf Sp}_{2l}}^1\rightarrow M_{{\bf SL}_{2l}}^1$.

\th Proposition
\enonce {\rm a)} The group $\Pic(M_{{\bf
Sp}_{2l}}^{1})$ is generated by ${\cal L}$.
\ind {\rm b)} The group $\Pic(M_{{\bf
PSp}_{2l}}^{1})$ is generated by ${\cal L}$ if $l$ is even, and by ${\cal
L}^2$ if $l$ is odd.
\endth
\label{C}
 {\it Proof}: By Prop.\ \ref{cyc} to prove a) it suffices to prove that
${\cal L}$ is not divisible. Choose an element $(A,\psi )$ of  $M_{{\bf
Sp}_{2l-2}}^{1}$, and consider the map $u:M_{{\bf SL}_2}^1\rightarrow M_{{\bf
Sp}_{2l}}^{1}$ given by $u(E)=(E,\det )\oplus (A,\psi )$. The pull back
of ${\cal L}$ is the theta line bundle $\Theta $ on $M_{{\bf SL}_2}^1$, hence
the
assertion a).  \ind Let us prove b).
By Remark \ref{M(2l,l)} we already know that ${\cal L}^2$ descends to  $M_{{\bf
PSp}_{2l}}^{1}$, and that ${\cal L}$ descends if $l$ is even.
 Consider the morphism $\mu
:M_{{\bf SL}_2}^1\rightarrow M_{{\bf Sp}_{2l}}^{1}$ given by $\mu (E
)=(E,\det)^{\oplus l}$. One has $\mu ^*{\cal L}=\Theta ^l$, so if ${\cal L}$
descends
 $\Theta ^l$ descends to $M_{{\bf PGL}_2}^1$; by Prop.\
\ref{A_d} this implies that $l$ is even.\cqfd

\medskip
\subsection Let us consider the group $G={\bf PSO}_{2l}$.
The moduli space $M_{G}$ has $4$ components, indexed by the center
$\{1,-1,\varepsilon
,-\varepsilon \}$ of ${\bf Spin}_{2l}$ (\ref{PSO}).
\ind The component $M^1_{{\bf PSO}_{2l}}$  has already been dealt with in
Prop.\
\ref{BC0}. The component $M_{{\bf PSO}_{2l}}^{-1}$
 is the quotient by the action of $J_2$ of the moduli space
$M^{-1}_{{\bf SO}_{2l}}$ of semi-stable quadratic bundles with $w_2=1$. Let
${\cal D}$ denote the determinant bundle on this moduli space.

\smallskip
\th Proposition
\enonce The group $\Pic(M_{{\bf PSO}_{2l}}^{-1})$ is generated by ${\cal
D}^2$ if  $l$ is even, by ${\cal D}^4$ if $l$ is odd.
\endth
{\it Proof}: The same proof as in \ref{BC0} shows that ${\cal D}^2$
descends to
$M_{{\bf PSO}_{2l}}^{-1}$ if $l$ is even, and that ${\cal D}^4$ descends
but ${\cal D}^2$ does not if $l$ is odd. To prove that ${\cal D}$ does not
descend when $l$ is even $\ge 3$, we  apply the argument of {\it loc.\ cit.}
to the morphism $u:JX\times M_{{\bf SO}_{2l-2}}^{-1}\rightarrow M_{{\bf
SO}_{2l}}^{-1}$ deduced from the natural embedding ${\bf SO}_2\times{\bf
SO}_{2l-2}\mono {\bf SO}_{2l}$ (note that $w_2$ is additive and $w_2(\alpha
\oplus\alpha ^{-1})=0$ for $\alpha \in JX$).  \ind When $l=2$ we consider
the morphism $v:M_{{\bf SL}_2}^1\times M_{{\bf SL}_2}^1\longrightarrow
M_{{\bf SO}_4}^{-1}$  which associates to a pair $(E,F)$  the vector bundle
${\cal H}om(E,F)$ with the quadratic form defined by the determinant and
the orientation deduced from the canonical isomorphism  $\det(E^*\otimes
F)\iso$ $(\det E)^{-2}\otimes(\det F)^2$. One has $v^*{\cal D}={\cal L}
\boxtimes {\cal L}$, where ${\cal L}$ is the theta line bundle on $M_{{\bf
SL}_2}^1$. Since ${\cal L}$ does not descend to $M^1_{{\bf PGL}_2}$ (Prop.\
\ref{A_d}), it follows from the commutative diagram $$\diagramme{M_{{\bf
SL}_2}^1\times M_{{\bf SL}_2}^1 &\hfl{v}{} &M_{{\bf SO}_4}^{-1} \cr
\vfl{}{} && \vfl{}{}\cr M^1_{{\bf PGL}_2}\times M^1_{{\bf PGL}_2}& \hfl{}{}
& M_{{\bf PSO}_4}^{-1} }$$that ${\cal D}$  does not descend to $M_{{\bf
PSO}_4}^{-1}$.\cqfd

\medskip
\ind We now consider
 the components  $M^{\pm\varepsilon }_{{\bf PSO}_{2l}}$ corresponding to the
 elements $+\varepsilon $ and $-\varepsilon $ of the center of ${\bf
Spin}_{2l}$. Each of these is the quotient by $J_2$ of the moduli space
$M^{\pm\varepsilon }_{{\bf SO}_{2l}}$
of semi-stable quadratic bundles $(E,q,\omega )$, where  $E$ is a vector bundle
of rank $2l$,  $q:\sym^2E\rightarrow {\cal O}_X(p)$  a quadratic form and
$\omega
:\det E\rightarrow {\cal O}_X(lp)$ an isomorphism compatible with  $q$;
changing the
sign of $\omega $ exchanges   $M^\varepsilon $ and
$M^{-\varepsilon }$ (\ref{PSO}).
 We  denote by ${\cal L}_l$ the pull back of the theta line bundle  on
$M_{{\bf
SL}_{2l}^l}$ under the natural map $M^{\pm\varepsilon }_{{\bf
SO}_{2l}}\rightarrow
M_{{\bf SL}_{2l}^l}$.

\th Proposition
\enonce The group $\Pic(M^{\pm\varepsilon }_{{\bf PSO}_{2l}})$ is generated by
${\cal
L}_l$ when $l$ is even, and by ${\cal L}_l^2$ when $l$ is odd. \endth
{\it Proof}: We already know that the theta line bundle descends to $M_{{\bf
SL}_{2l}}^l/J_2$ when $l$ is even  and that its square
descends when $l$ is odd (\ref{M(2l,l)}), so we have only to prove that ${\cal
L}_l$
does not descend when $l$ is odd.  \ind Let us first consider the case
$l=3$. If $E$ is a vector bundle of rank $4$ and determinant ${\cal O}_X(p)$ on
$X$,
the bundle $\wedge^2E$  carries a quadratic form with values in ${\cal O}_X(p)$
(defined by the exterior product) and an orientation. We thus get a morphism
 $\lambda : M_{{\bf SL}_4}^1\rightarrow M^{\pm\varepsilon }_{{\bf SO}_6}$ such
that
$\lambda (E\otimes\alpha )=\lambda (E)\otimes\alpha ^2$ for $\alpha \in J_4$.
 An easy computation  shows that $\lambda ^*{\cal L}_3$ is the theta line
bundle
on $M_{{\bf SL}_4}^1$, which does not descend to $M_{{\bf PGL}_4}^1$
(\ref{A_d});  our assertion follows. \ind For $l$ odd $\ge 5$, we consider the
morphism $\mu :M^{\pm\varepsilon }_{{\bf SO}_{2l-6}}\times M^{\pm\varepsilon
}_{{\bf SO}_6}\longrightarrow M^{\pm\varepsilon }_{{\bf SO}_{2l}}$ deduced from
the embedding ${\bf SO}_{2l-6}\times {\bf SO}_6\mono  {\bf SO}_{2l}$. It is
$J_2$\tx equivariant (with respect to the canonical action of $J_2$ on the
spaces $M^{\pm\varepsilon }_{{\bf SO}_{2n}}$, and the diagonal action on the
product), and the pull back $\mu ^*{\cal L}_l$ is isomorphic to ${\cal
L}_{l-3}\boxtimes {\cal L}_3$.  Assume that ${\cal L}_l$ descends to
$M^{\pm\varepsilon }_{{\bf PSO}_{2l}}$; since ${\cal L}_{l-3}$ descends, we
deduce from \ref{produit} that ${\cal L}_3$ descends, contradicting what we
just
proved.\cqfd

\vskip 1cm

\section {Determinantal line bundles}
\subsection \label{e_G} We can express the above results in a more suggestive
way.
Assume $G$ is of type $A,B,C$ or $D$; let $\delta\in \pi_1(G)$. We identify
$\Pic(M_G^\delta)$ to a subgroup of $\Pic({\cal
M}_{\widetilde{G}}^{\delta,ss})$. Let
$\sigma $ be the standard representation of $\widetilde{G}$ in ${\bf C}^r$ (for
$\widetilde{G}={\bf SL}_r$, ${\bf Spin}_r$ or ${\bf Sp}_r$), and ${\cal
D}_\sigma $
the corresponding determinant bundle on ${\cal M}_{\widetilde{G}}^{\delta,ss}$.
The
results of sections 8 to 11 express the generator of $\Pic(M_G^\delta)$ as a
certain
power of
 ${\cal D}_\sigma$. Using the fact that the  pull back to ${\cal
M}^{d,ss}_{{\bf SL}_r}$ of the theta line  bundle on $M_{{\bf SL}_r }^d$  is
$({\cal D}_\sigma )^{r\over (d,r)}$, one finds:
\th Proposition
\enonce Assume that $G$ is one of the groups ${\bf PGL}_r$, ${\bf PSp}_{2l}$ or
${\bf
PSO}_{2l}$. Put
$\varepsilon_G^{}=1$ if  the rank of $G$ is even, $2$ if it is odd. Let
$\delta\in \pi_1(G)$. The group $\Pic(M_G^\delta)$ is generated by
$({\cal D}_\sigma )^{r\varepsilon_G^{}}$ for $G={\bf PGL}_r$, and by $({\cal
D}_\sigma)^{2\varepsilon_G^{}}$
 for the other groups.
\endth

\medskip

 \subsection  To produce line bundles on
$M_G^\delta$, we have already used  the following recipe: to
any representation
 $\rho :G\rightarrow {\bf SL}_N$ we associate the pull back ${\cal D}_\rho $ of
the
determinant bundle under the morphism  $M_G^\delta \rightarrow M_{{\bf SL}_N}$
deduced from $\rho $. These line bundles generate a subgroup $\Pic_{\rm
det}(M_G^\delta)$ of $\Pic(M_G^\delta)$. We suspect that this subgroup is
actually
equal to $\Pic(M_G^\delta)$, i.e.\ that all line bundles on $M_G^\delta$ can be
constructed from representations of $G$. We have checked this in some cases:

\th Proposition
\enonce Assume $G$ is of classical type or of type $G_2$, and either simply
connected or adjoint or isomorphic to ${\bf SO}_r$. Then, for every $\delta\in
\pi_1(G)$, one has $\Pic_{\rm det}(M_G^\delta)=\Pic(M_G^\delta)$.
\endth
{\it Proof}: The simply connected case, and also the case $G={\bf SO}_r$,
follow
from [L-S], Prop.\ \ref{Pic_Spin} and \ref{G_2}.
\ind  The other groups are those which appear in the above Proposition; let
us denote by $e_G$ the positive integer such that $({\cal D}_\sigma)^{e_G}$
generates $\Pic(M_G^\delta)$. If
$\rho
$ is a representation of $G$, with Dynkin index $d_\rho $, the line bundle
${\cal D}_\rho $ on $M_{G}^d$ is isomorphic to $({\cal D}_\sigma) ^{
d_\rho/d_\sigma }$ ($d_\sigma $ is $1$ for the types $A,C$ and $2$ for
$B,D$). It follows that $e_G$ divides $d_\rho /d_\sigma$, and that our
assertion is equivalent to saying that $e_G$ is the g.c.d.\ of the numbers
$d_\rho/d_\sigma$ when
 $\rho$ runs over the representations of $G$.

\ind Let us consider the case $G={\bf
PGL}_r$.  We have $d_{\rm Ad}=2r$, which settles the case $r$ even.
If $r$ is odd, consider the representation $\sym^2 \otimes \wedge^{r-2}$
of ${\bf SL}_r$; since $\Mu_r$ acts trivially, it defines a representation
$\rho $ of
${\bf PGL}_r$, whose Dynkin index is
$$\eqalign{d_\rho &=d_{\sym^2}\,\dim
\wedge^{r-2}+d_{\wedge^2}\,\dim \sym^2 \cr
&=(r+2){r\choose 2}+(r-2){r+1\choose
2} = r^3-2r\ .\cr}$$
Then $(d_{\rm Ad},d_\rho )=r=e_G$, which proves the result in this case.

\ind For $G={\bf PSp}_{2l}$, easy computations give
$d_{\rm Ad}=2l+2$ and $d_{\wedge^2}=2l-2$, hence $e_G=(d_{\rm
Ad},d_{\wedge^2})$. For
$G={\bf PSO}_{2l}$, one has
$d_{\rm Ad}=2(2l-2)$ and $d_{\sym^2}=2(2l+2)$, hence $e_G=(d_{\rm
Ad},d_{\wedge^2})/d_\sigma $.\cqfd

\medskip
\rem{Remark} We can also prove the equality $\Pic_{\rm
det}(M_G^0)=\Pic(M_G^0)$ for $G={\bf SL}_r/\Mu_s$ when $s$ and $r/s$ are
coprime. Reasoning as above and using Prop. \ref{A0}, we need to prove that the
g.c.d.\  of the $d_\rho $'s is $2s$ if $s$ is even, and $s$ if it is odd.
 We consider the representation $\rho _p=\sym^p\otimes \wedge^{s-p}$ for
$1\le p \le s$. Using some nontrivial combinatorics we can prove the relation
$\displaystyle \sum_{p=1}^s p\,d_{\rho _p}=(-1)^ss^2$. Since $d_{\rm Ad}=2r$
we find  that the g.c.d.\  of the $d_\rho $' divides $(2r,s^2)=s(2{r\over
s},s)$, hence our assertion.\cqfd
\vskip1cm
\section{Local properties of the moduli spaces $M_G$}
\def\Cl{\mathop{\rm Cl}\nolimits}

\subsection  A  $G$\tx bundle $P $ is called {\it regularly
stable} if it is stable and its automorphism group is equal to the center
$Z(G)$ of $G$. The open subset $M_G^{\rm reg}$ of $M_G$ corresponding to
regularly stable $G$\tx bundles is smooth, and its complement in $M_G$ is of
codimension $\ge 2$, except when $X$ is of genus $2$ and $G$ maps onto ${\bf
PGL}_2$: this is seen exactly as the analogous statement for Higgs bundles,
which is proved in [F1], thm.\ II.6. In what follows we will assume that we
are not in this exceptional case, leaving to the reader to check that our
assertions extend by using the
explicit  description of $M_{{\bf SL}_2}$ in genus $2$.
\ind Let
$i$ be the natural injection of $M_G^{\rm reg}$ into $M_G$. Then the
 map $i_*$ identifies $\Pic(M_G^{\rm reg})$ with the Weil divisor class
group $\Cl(M_G)$, that is the group  of isomorphism classes of rank $1$
reflexive sheaves on
$M_G$ (see [Re], App.\ to \S 1); the restriction map
$i^*:\Pic(M_G)\rightarrow \Pic(M_G^{\rm reg})$ corresponds to the inclusion
$\Pic(M_G) \i \Cl(M_G)$. Local factoriality of M$_G$ is equivalent to the
equality $\Pic(M_G) = \Cl(M_G)$.

\ind   We already know from [D-N] and [L-S] that $M_G$ is locally
factorial when $G$ is ${\bf SL}_r$ or ${\bf Sp}_{2l}$. We want to show that
these are essentially the only cases where this occurs.

\th Proposition
\enonce Let $G$ be a simply connected group, containing a factor of type
$B_l\ (l\ge 3)$, $D_l\ (l\ge 4)$, $F_4$ or $G_2$. Then $M_G$ is not locally
factorial. \endth
\ind We believe that this still holds if $G$ contains a factor of type
$E_l$. This would have the amusing consequence that {\it the semi-simple
groups $G$ for which $M_G$ is locally factorial are exactly those which
are} special {\it in the sense of Serre, i.e.\ such that all  $G$\tx
bundles are locally trivial for the Zariski topology} (see [G2]). \par

\smallskip
{\it Proof of the Proposition}: We can assume that $G$ is almost simple.
Choose  a presentation of
${\cal M}_G$ as a quotient of a smooth scheme $R$ by a reductive group
$\Gamma $, such that
$M_G$ is a good  quotient of $R$ by $\Gamma $   (lemma \ref{comppres}).
We denote by $\sigma $ the
 standard representation
in ${\bf C}^r$ in case $G={\bf Spin}_r$,  in ${\bf
C}^7$ if $G$ is of type $G_2$, and the orthogonal representation in ${\bf
C}^{26}$ with highest weight $\varpi_4$ if $G$ is of type $F_4$ (we use the
standard notation of [Bo], Lie VII). Let ${\cal D}$ be the determinant
bundle on $R$ associated to $\sigma $. As in the proof of Prop.\
\ref{Pic_Spin},  the choice of a theta-characteristic $\kappa $ on $X$
allows us to define a square root ${\cal P}_\kappa $ of
${\cal D}$ on $R$, with a canonical $\Gamma $\tx linearization. We will
show that ${\cal P}_\kappa $ descends to the open subset $M_G^{\rm reg}$, but
not
to $M_G$, thus showing that the restriction map is not surjective.
\ind The
first assertion is clear if $G$ is of type $F_4$ or $G_2$, because then
$Z(G)$  is trivial, so $\Gamma $ acts freely on the open subset of $R$
corresponding to regularly stable $G$\tx bundles. Suppose $G={\bf Spin}_r$;
let $Q$ be a $G$\tx bundle, and $z$ an element of $Z(G)$. The image of $z$
in ${\bf SO}_r$ is either $1$ or possibly $-1$ if $r$ is even;   since
$h^0(Q_\sigma \otimes\kappa )\equiv rh^0(\kappa )$ (mod.$\,2$) by [L-S],
7.10.1, we conclude that $z$ acts trivially on $\wedge^{\rm max}
H^0(Q_\sigma \otimes\kappa )$, i.e.\ on the fibre of ${\cal P}_\kappa $ at
$Q$ (\ref{Pic_Spin}).

\ind We already know that ${\cal P}_\kappa $ does not descend to $M_G$ when
$G={\bf Spin}_r$ (Prop.\ \ref{Pic_Spin}) or $G$ is of type  $G_2$ (Prop.\
\ref{G_2}); it remains to show that the class of ${\cal D}$ is not
divisible by $2$ in $\Pic(M_G)$ when $G$ is of type $F_4$. There is a
natural inclusion ${\bf Spin}_8\i G$, which induces a morphism $f:M_{{\bf
Spin}_8}\rightarrow M_G$.  An easy computation gives that the Dynkin index
of the restriction to ${\bf Spin}_8$ of the standard representation of $G$
is $6$. Since the Dynkin index of the standard representation of ${\bf
Spin}_8$ is $2$, it follows that $f^*{\cal D}$ is  isomorphic to ${\cal
D}_0^{\otimes 3}$, where ${\cal D}_0$ is the  generator $\Pic(M_{{\bf
Spin}_8})$; this implies that ${\cal D}$ is not divisible by $2$ in
$\Pic(M_G)$.\cqfd
 \medskip

\ind We now treat the case of a non simply connected group. We start with
two  lemmas which are certainly well known, but for which we could find no
reference:

\th Lemma
\enonce Let $\pi :\widetilde{Y}\rightarrow Y$ be a ramified  Galois covering,
with abelian Galois group $A$. If $\pi $ is \'etale in codimension $1$, the
variety $Y$ is not locally factorial.
\endth\label{factorial}
{\it Proof}: Let $Y^{\rm o}$ be an open subset of $Y$ such that $Y\moins
Y^{\rm o}$ has codimension $\ge 2$ and the induced covering
$\pi ^{\rm o}:\widetilde{Y}^{\rm o}\rightarrow Y^{\rm o}$ is \'etale. This
covering corresponds to a homomorphism $L: \widehat{A}\rightarrow
\Pic(Y^{\rm o})$ such that $\pi _*{\cal O}_{\widetilde{Y}^{\rm o}}=\sdir_{\chi
\in
\widehat{A}}^{}L(\chi) $. If $Y$ is locally factorial, the
restriction map $\Pic(Y)\rightarrow \Pic(Y^{\rm o})$ is bijective, so $L$
extends
to a  homomorphism $ \widehat{A}\rightarrow \Pic(Y)$ which defines an \'etale
covering of $Y$ extending $\pi ^{\rm o}$, and therefore equal to $\pi $.
Then $\pi $ is \'etale, contrary to our hypothesis.\cqfd

\th Lemma
\enonce Let $S$ be a scheme, $H$ an algebraic group, $A$ a closed central
subgroup of $H$, $P$ a principal $H$\tx bundle on $S$. The cokernel of the
natural homomorphism $\Aut(P)\rightarrow \Aut(P/A)$ is canonically isomorphic
to
the stabilizer of $P$ in $H^1(X,A)$ {\rm (}for the natural action of $H^1(X,A)$
on
$H^1(X,H)\,)$.
\endth
\label{aut}
{\it Proof}: Denote by ${\cal A}ut\,(P)$ the automorphism  bundle  of the
$H$\tx
bundle $P$. We have an exact sequence of groups over $S$
$$1\rightarrow A_S \longrightarrow {\cal A}ut\,(P)\longrightarrow {\cal
A}ut\,(P/A)\rightarrow 1$$(to check exactness one may replace $P$ by the
trivial $H$\tx bundle, for which this is clear). The associated cohomology
exact sequence reads
$$1\rightarrow A \longrightarrow \Aut(P)\longrightarrow \Aut(P/A)
\longrightarrow H^1(S,A)\qfl{h} H^1(S,{\cal A}ut\,(P))\ .$$
The map $h$ associates to an $A$\tx bundle $\alpha$ the class of the ${\cal
A}ut\,(P)$\tx bundle
$\alpha\times^A{\cal A}ut\,(P)$, which is canonically isomorphic to ${\cal
I}som\,(P,\alpha\times^AP)$; the element $h(\alpha)$ is trivial if and only if
this last bundle admits a global section, which means exactly that
$\alpha\times^AP$ is isomorphic to $P$, hence the lemma.\cqfd

\th Proposition \enonce Suppose $G$ is not simpy
connected; let $\delta\in \pi _1(G)$. The moduli space
$M_G^\delta$ is not locally factorial.
\endth
{\it Proof}: We first prove that the Galois covering $\pi
:M_{\widetilde{G}}^\delta\rightarrow M_G^\delta$ is \'etale above
$(M_G^\delta)^{\rm reg}$. We put $A=\pi _1(G)$,  and choose an isomorphism
$A\iso
\pprod_{j=1}^s\Mu_{r_j}$; we use freely the notation of
(\ref{M_G^\delta}). We denote by
$H$ the group
$C_A\widetilde{G}=(\widetilde{G}\times T)/A$.
 Let $Q\in (M_G^\delta)^{\rm reg}$ and  $P$ a point of
$M_{\widetilde{G}}^\delta$ above $Q$; we will use the same letters to denote
the corresponding bundles. Using the isomorphism $H/A\cong G\times (T/A)$,
the condition $\pi(P)=Q$ means that the
$(H/A)$\tx bundle $P/A$ is isomorphic to $Q\times {\cal O}_X({\bf d}p)$.
Since $\Aut(Q)$ is reduced to the
center of  $G$,  the map $\Aut(P)\rightarrow
\Aut(P/A)$ is surjective; we deduce from lemma  \ref{aut} that the
stabilizer of $P$ in $H^1(X,A)$ is trivial, i.e.\ $\pi$ is \'etale at
$P$.
\ind It follows that the abelian cover $\pi:M_{\widetilde{G}}^{\delta}\ra
M_{G}^{\delta}$ is \'etale in codimension one. Since it is
ramified by lemma \ref{unirat}, we conclude from lemma
\ref{factorial} that $M_{G}^{\delta}$ is not locally
factorial.

\medskip
\ind Finally we observe that, though the moduli space is not locally
factorial in most cases, it is always Gorenstein (this is proved in [K-N],
thm.\ 2.8, for a simply connected $G$):
\th Proposition
\enonce The moduli space $M_G$ is Gorenstein.
\endth
{\it Proof}: We choose again a presentation of ${\cal M}_G$ as a quotient
of a smooth scheme $R$ by a reductive group $\Gamma $, such that $M_G$ is a
good quotient of $R$ by $\Gamma $  (lemma \ref{comppres}); we denote by
${\cal P}$ the universal bundle on $X\times R$, and by  $R^{\rm reg}$
the open subset of $R$ corresponding to regularly stable bundles.
 Since the center of $G$ is killed by the adjoint representation,
 the vector bundle
$\ad({\cal P})$ descends to a vector bundle on $X\times M_G^{\rm reg}$,
that we will still denote $\ad({\cal P})$.
Deformation theory provides an isomorphism
$T_{M_G^{\rm reg}} \iso R^1pr_{2*}(\ad {\cal
P})$; since $H^0(X,\ad(P))=0$ for $P\in M_G^{\rm reg}$,
the line bundle $\det T_{M_G^{\rm reg}}$ is isomorphic to $\det
Rpr_{2*}(\ad {\cal
P}))$, that is to the restriction to $M_G^{\rm reg}$ of the determinant
bundle  ${\cal D}_{\rm Ad}$
associated to the adjoint representation.
\ind  Since $M_G$ is
Cohen-Macaulay, it admits a dualizing sheaf $\omega $, which is
torsion-free and reflexive ([Re], App. of \S 1). The reflexive sheaves
$\omega $ and ${\cal D}_{\rm Ad}^{-1}$, which are isomorphic above
$M_G^{\rm reg}$, are isomorphic ({\it loc. cit.}), hence $\omega $ is
invertible.\cqfd

\vskip2cm \centerline{ REFERENCES} \bigskip \baselineskip13pt
\def\num#1{\item{\hbox to\parindent{\enskip [#1]\hfill}}}
\parindent=1.5cm

\num{B-L} A.\ {\pc BEAUVILLE}, Y.\ {\pc LASZLO}: {\sl
Conformal blocks and generalized theta functions.} Comm.
Math.\ Phys.\  {\bf 164}, 385-419 (1994).
\smallskip
\num{B-N-R} A.\ {\pc BEAUVILLE}, M.S.\ {\pc NARASIMHAN}, S.\ {\pc RAMANAN}:
{\sl  Spectral curves and the generalised theta divisor}. J.\ reine angew.\
Math.\ {\bf 398}, 169-179 (1989). \smallskip
\num{Bo} N.\ {\pc BOURBAKI}: {\sl El\'ements de Math\'ematique}. Hermann,
Paris.   \smallskip

\num{Br} K.\ {\pc BROWN}: {\sl Cohomology of groups}. GTM {\bf
87}, Springer-Verlag (1982).
\smallskip
 \num{D-N} J.M.\ {\pc DREZET}, M.S.\ {\pc NARASIMHAN}: {\sl Groupe de Picard
des vari\'et\'es de modules de fibr\'es semi-stables sur les courbes
alg\'ebriques.} Invent.\ math.\ {\bf 97}, 53-94 (1989).
 \smallskip

\num {F 1} G.\ {\pc FALTINGS}: {\sl Stable $G$\tx bundles and projective
connections.} J.\ Algebraic Geometry {\bf 2}, 507-568 (1993). \smallskip
\num {F 2} G.\ {\pc FALTINGS}: {\sl A proof for the Verlinde formula.} J.\
Algebraic Geometry {\bf 3}, 347-374 (1994). \smallskip

\num {G 1} A.\ {\pc GROTHENDIECK}: {\sl Techniques de construction et
d'existence en g\'eom\'etrie
alg\'ebrique: les sch\'emas de Hilbert.}
S\'em.\ Bourbaki {\bf 221},  1-28 (1960/61).
\smallskip

\num {G 2} A.\ {\pc GROTHENDIECK}: {\sl Torsion homologique et sections
rationnelles.}
S\'emi\-naire Chevalley, $2^e$ ann\'ee, Exp.\ 5. Institut Henri Poincar\'e,
Paris
(1958). \smallskip
\num{Ha} A.\ {\pc HARDER}: {\sl Halbeinfache Gruppenschemata \"uber
Dedekindringen}. Invent.\ Math.\ {\bf 4}, 165-191 (1967).
\smallskip

\num{H} H.\ {\pc HIRONAKA}: {\sl Introduction to the theory of infinitely near
singular points}. Memorias de Mat.\ del Inst.\ ``Jorge Juan", Madrid
(1974).\smallskip

\num{Hi} A.\ {\pc HIRSCHOWITZ}: {\sl Probl\`emes de
Brill-Noether en rang sup\'erieur}. C.\ R.\ Acad.\ Sci.\ Paris {\bf
307}, 153-156 (1988).
\smallskip
 \num{K-N} S.\ {\pc KUMAR}, M.S.\ {\pc NARASIMHAN}: {\sl
Picard group of the moduli spaces of $G$\tx bundles}.
Preprint alg-geom/9511012.\smallskip

\num{K-N-R} S.\ {\pc KUMAR}, M.S.\ {\pc NARASIMHAN}, A.\ {\pc RAMANATHAN}:
{\sl Infinite Grassmannians and moduli spaces of $G$\tx
bundles}. Math.\ Annalen {\bf 300}, 41-75 (1994).
\smallskip \num{L} H.
{\pc LANGE}: {\sl Zur Klassification von Regelmannigfaltigkeiten}.
Math.\ Annalen {\bf 262}, 447-459 (1983). \smallskip

\num{L-MB} G.\ {\pc LAUMON}, L.\ {\pc MORET-BAILLY}:  {\sl Champs
alg\'ebriques.} Pr\'epublication Universit\'e Paris-Sud (1992).
\smallskip

\num{L-S} Y.\ {\pc LASZLO}, Ch.\ {\pc SORGER}: {\sl The line bundles
on the moduli of parabolic $G$\tx bundles over curves and their
sections}. Annales de l'ENS, to appear; preprint alg-geom/9507002.
\smallskip

\num{M} O.\ {\pc MATHIEU}: {\sl Formules de caract\`eres pour les alg\`ebres
de Kac-Moody g\'en\'e\-rales.} Ast\'erisque {\bf 159-160} (1988).
\smallskip

\num{N-R} M.S.\ {\pc NARASIMHAN}, S.\ {\pc RAMANAN}:
{\sl Generalized Prym varieties as fixed points}. J.\ of the Indian
Math.\ Soc.\ {\bf 39}, 1-19 (1975).
\smallskip
\num{O} W.M.\ {\pc OXBURY} : {\sl Prym varieties and the moduli of spin
bundles}. Preprint, Durham
(1995).

\smallskip
\num{R1} A.\ {\pc RAMANATHAN}: {\sl Stable principal $G$\tx bundles}. Thesis,
Bombay University (1976).
\smallskip

\num{R2} A.\ {\pc RAMANATHAN}: {\sl Stable principal bundles on a compact
Riemann surface}. Math.\ Ann.\ {\bf 213}, 129-152 (1975).
\smallskip
\num{Re} M.\ {\pc REID}: {\sl Canonical $3$-folds}. Journ\'ees de
G\'eom\'etrie alg\'ebrique d'Angers (A. Beauville ed.), 273-310; Sijthoff and
Noordhoff (1980).
\smallskip
\num{Se} C.S.\ {\pc SESHADRI}: {\sl Quotient spaces modulo reductive algebraic
groups}. Ann.\ of Math.\ {\bf 95}, 511-556 (1972).
\smallskip

 \num{S1} R.\ {\pc STEINBERG}: {\sl G\'en\'erateurs, relations et
rev\^etements de groupes alg\'e\-briques}. Colloque sur la th\'eorie des
groupes
alg\'ebriques (Bruxelles), 113-127; Gauthier-Villars, Paris (1962).
\smallskip
\num{S2} R.\ {\pc STEINBERG}:
{\sl Regular elements of semisimple algebraic groups}. Publ.\ Math.\ IHES {\bf
25}, 281-312 (1965).
\smallskip
\def\pc#1{\eightrm#1\sixrm}
$$\hbox to 16truecm{\eightrm\vtop{\hbox to 5cm{\hfill A. {\pc BEAUVILLE}, Y.
{\pc LASZLO}\hfill}
 \hbox to 5cm{\hfill DMI -- \'Ecole Normale
Sup\'erieure\hfill} \hbox to 5cm{\hfill (URA 759 du CNRS)\hfill}
\hbox to 5cm{\hfill  45 rue d'Ulm\hfill}
\hbox to 5cm{\hfill F-75230 {\pc PARIS} Cedex 05\hfill}}
\vtop{\hbox to 5cm{\hfill Ch. {\pc SORGER}\hfill}\par
\hbox to 5cm{\hfill Institut de Math\'ematiques de Jussieu\hfill}\par
\hbox to 5cm{\hfill (UMR 9994 du CNRS)\hfill}
\hbox to 5cm{\hfill Univ. Paris 7 -- Case Postale 7012\hfill}
\hbox to 5cm{\hfill 2 place Jussieu\hfill}
\hbox to 5cm{\hfill F-75251 {\pc PARIS} Cedex 05\hfill}}} $$

\end